\documentclass[11pt,
]{article}

\usepackage[tbtags]{amsmath}
\usepackage{amsfonts,amssymb,mathrsfs,amscd,amstext,
latexsym,amsbsy,amsthm}

\usepackage{bbold}
\usepackage{mparhack}
\usepackage{makeidx}
\usepackage{marginnote}

\input{34.sty}

\makeindex

\begin{document}

\title
{Definable $\Eo$ classes at arbitrary projective levels
}

\author 
{Vladimir Kanovei\thanks{IITP RAS and MIIT,
 \ {\tt kanovei@googlemail.com} --- 
contact author. 
Partial support of grant RFBR 17-01-00705 acknowledged.
}  
\and
Vassily Lyubetsky\thanks{IITP RAS,
\ {\tt lyubetsk@iitp.ru}. 
Partial support of grant RSF 14-50-00150 acknowledged. 
}}

\date 
{\today}

\maketitle


\begin{abstract}

Using a modification of the invariant Jensen forcing 
of \cite{kl22}, 
we define a model of $\zfc$, in which, for a given $n\ge3$, 
there exists a lightface 
\dd{\ip1n}set of reals, 
which is a \dd\Eo equivalence class, hence a 
countable set, and which does  
not contain any $\od$ element, while 
every non-empty countable \dd{\is1{n}}set of reals 
is necessarily constructible, hence contains only 
$\od$ reals. 
\end{abstract}

\markright{Definable $\Eo$ classes at arbitrary projective levels}

{\small
\def\contentsname{\large Contents}
\tableofcontents
}

\parf{Introduction}
\las{or}

Problems related to definability of mathematical objects, 
were one of focal points of the famous discussion on 
mathematical foundations in the beginning of XIX C. 
In particular, Baire, Borel, Hadamard, and Lebesgue, 
participants of the exchange of letters published in 
\cite{cinq}, in spite of essential disagreement between 
them on questions related to mathematical foundations, 
generally agreed that the proof of existence of 
an element in a given set, and a direct definition  
(or effective construction) of such an element  
--- are different mathematical results, of which  
the second does not follow from the first. 
In particular, Lebesgue in his contribution to 
\cite{cinq} pointed out the difficulties in the problem  
of effective choice, that is, choice of definable element 
in a definable non-empty set.%
\snos{Ainsi je vois d\'ej\`a une difficult\'e dans ceci 
\lap{dans un $M'$ d\'etermin\'e
je puis choisir un $m'$ d\'etermin\'e}, in the French 
original. 
Thus I already see a difficulty with the assertion that 
\lap{in a determinate $M'$ I can choose a determinate $m'$}, 
in the translation.}

Studies in modern set theory demonstrated that effective 
effective choice is not always possible. 
In particular, it is true in many well-known models  
(including the very first Cohen models),  
that the set $X=\dR\bez\rL$ of all 
Goedel-nonconstructible reals is not empty, but  
contains no definable elements. 

One may note that if the set $X$ is non-empty then it 
has to be rather large, that is, surely of cardinality 
$\mathfrak c$, if measurable then of full measure, \etc. 
Is there such an example among \rit{small}, \eg\ 
countable sets? 
This problem was discussed at 
\rit{Math\-over\-flow}\snos
{\label{snos1}
A question about ordinal definable real numbers. 
\rit{Math\-overflow}, March 09, 2010. 
{\tt http://mathoverflow.net/questions/17608}. 
}
and 
\rit{Foundations of mathematics} (FOM)\snos
{\label{snos2}%
Ali Enayat. Ordinal definable numbers. FOM Jul 23, 2010.
{\tt http://cs.nyu.edu/pipermail/fom/2010-July/014944.html}}
. 

The problem was solved in \cite{kl27a} 
(to appear in \cite{kl27}). 
Namely, let $\rL[\sis{a_n}{n<\om}]$ be a 
\dd{\jf^\om}generic extension of $\rL$, where 
$\jf^\om$ is the countable power 
(with finite support)
of \rit{Jensen's minimal forcing} $\jf$ \cite{jenmin}%
\index{forcing!Jensen's forcing, $\jf$}%
\index{zzJ@$\jf$}%
\snos
{See also 28A в \cite{jechmill} on this forcing. 
We acknowledge that the idea to use the countable
power $\jf^\om$ of Jensen's forcing $\jf$ to
define such a model belongs to
Ali Enayat~\cite{ena}.}.
The key property of $\jf$ is that it adds 
a nonconstructible generic $\id13$ real $a\in\dn$, 
in fact $\ans a$ is a $\ip12$ singleton.
Accordingly $\jf^\om$ adds a sequence $\sis{a_n}{n<\om}$ 
of \dd\jf generic reals real $a_n\in\dn$ to $\rL$, and 
as shown in \cite{kl27a,kl27} there is no orher  
\dd\jf generic reals in $\rL[\sis{a_n}{n<\om}]$ except 
for the reals $a_n$. 
Furthermore, since \lap{being a \dd\jf generic real} is 
a $\ip12$ relation, it is true in $\rL[\sis{a_n}{n<\om}]$ 
that $A=\ens{a_n}{n<\om}$ is a countable (infinite) 
lightface $\ip12$ set 
without $\od$ (ordinal-definable) elements.
\vyk{
\snos
{\od\ (ordinal-definable) is the class of all sets 
definable by formulas with ordinals as parameters
of the definition. 
This is the largest class of sets which can be called  
\rit{effectively definable} in the most general sense.} 
}

Using an uncountable product of forcing notions  
similar to $\jf^\om$,
we defined in \cite{kl28} a model in which 
there is a \lap{planar} $\ip12$ set with countable
vertical cross-sections, which cannot be uniformized  
by any real-ordinal definable (\rod) 
set.

For a more detailed analysis of the problem,
note that the elements $a_n$ of the set $A$, adjoined by  
the forcing $\jf^\om$, are connected to each other only by
the common property of their \dd\jf genericity. 
Does there exist a similar countable set with a more
definite mathematical structure? 

This question was answered in \cite{kl22}  
by a model in which there is 
an equivalence class of the equivalence
\index{equivalence relation!E0@$\Eo$}%
\index{zzE0@$\Eo$}%
relation $\Eo$\snos
{\label{deo}%
Recall that $\Eo$ is defined on the Cantor space $\dn$ 
so that $x\Eo y$ iff the set $\ens{n}{x(n)\ne y(n)}$
is finite. 
\dd\Eo equivalence classes are countable sets in $\dn$, 
of course.} 
(a \dd\Eo\rit{class}, for brevity), 
which is a (countable) lightface $\ip12$ set in $\dn$,
and does not contain  $\od$ elements.
This model makes use of a forcing notion $\jfi$, 
similar to Jensen's forcing $\jf$, but different from   
\index{forcing!Jensen's forcing, $\jf$!invariant, $\jfi$}%
\index{zzJi@$\jfi$}%
$\jf$. 
In particular, it consists of Silver trees  
(rather than perfect trees of general form, as $\jf$ does) 
and is invariant under finite transformations. 
Thus it can be called     
\rit{an invariant Jensen forcing}. 
Due to the invariance, $\jfi$ adjoins a 
\dd\Eo equivalence class of \dd{\jfi}generic reals 
rather than a single real, and this class turns out to 
be a $\ip12$ set without $\od$ elements.   
A forcing similar to $\jfi$ was also used in \cite{kl25}  
to define a model with a $\ip12$ Groszek -- Laver pair 
of \dd\Eo classes. 

Our main theorem extends this research line.  

\bte
\lam{Tun}
Let\/ $\nn\ge3$. 
There is a generic extension\/ $\rL[a]$ of\/ $\rL$, 
by a real\/ $a\in\dn$, such that the following is 
true in\/ $\rL[a]$$:$ 
\ben
\renu
\itlb{Tun1}%
$a\nin\od$ and\/ $a$ is minimal over\/ $\rL$ ---
hence each\/ $\od$ real belongs 
to\/ $\rL\,;$

\itlb{Tun2}%
the\/ \dd\Eo class\/ $\eko{a}$ is a 
(countable) 
lightface\/ $\ip1\nn$ set --- which by\/ \ref{Tun1} 
does not contain\/ $\od$ elements$;$

\itlb{Tun3}%
every countable\/ $\is1{\nn}$ set belongs to\/ $\rL$, 
hence consists of\/ $\od$ elements.
\een
\ete

\vyk{
\bte
\lam{Tun}
Let\/ $\nn\ge3$. 
There is a generic extension of\/ $\rL$, 
in which 
\ben
\renu
\itlb{Tun1}%
there exists a\/ $\ip1\nn$ 
\dd\Eo class not containing\/ $\od$ elements$;$

\itlb{Tun2}%
every countable\/ $\is1{\nn}$ set belongs to\/ $\rL$, 
hence consists of\/ $\od$ elements.
\een
\ete
}

Thus we have a model, in which, for a given $\nn\ge3$, 
there is a $\ip1\nn$, hence, $\od$ \dd\Eo class 
(therefore, a countable set) in $\dn$ 
containing only non-\od\ elements, and in the same time 
every countable $\is1{\nn}$ set belongs to\/ $\rL$, 
hence consists of $\od$ elements.
The abovementioned result of \cite{kl22} corresponds to
$\nn=2$ in this theorem since countable $\is12$ sets 
belong to $\rL$ and consist of $\od$ elements.

\parf{Connections to the Vitali equivalence relation}
\las{2vit}

The relation $\Eo$ is tightly connected with the
Vitali equivalence relation $\vit=\dR/\dQ$
\index{equivalence relation!Vitali@Vitali, $\vit$}%
\index{zzvit@$\vit$}%
on the true real line $\dR$.\snos
{$\vit$ is defined on $\dR$ so that
$x\vit y$ iff the set $x-y$
is rational.}
In particular, there is a lightface $\id11$
(in fact of a very low class $\id0n$)
injection $\vt:\rR\to\dn$ which \rit{reduces}
$\vit$ to $\Eo$ in the sense that the equivalence
${x\vit y}\eqv {\vt(x)\Eo\vt(y)}$ holds for all
$x,y\in \dR$. 
(See Mycielski and Osofsky \cite{myc} for an explicit
construction of $\vt$, along with an inverse
reduction of $\Eo$ to $\vit$.)
It follows that Theorem~\ref{Tun} is true with 
\ref{Tun2} for $[a]_{\vit}$
as well, since if $C\sq\dn$ is a $\ip1\nn$ 
\dd\Eo class not containing $\od$ elements then the
preimage $C'=\ens{x\in\dR}{\vt(x)\in C}$ is a $\ip1\nn$ 
Vitali class not containing $\od$ elements.

The interest to Vitali classes in this context is inspired
by the observation that they can be viewed as the most
elementary countable sets in $\ves$ which do not allow 
immediate effective choice of an element.
Indeed if a set $X\sq \ves$ contains at least one 
isolated, or even one-sided isolated point, then 
one of such points can be chosen effectively.
However any non-empty set $X\sq \ves$ without
one-sided isolated points is just an everywhere dense
set
(not counting close segments of the complementary set).
Yet the Vitali classes, that is, shifts of the rationals
$\raz$, are exactly the most simple and typical 
countable dense sets in $\ves$. 

Historically, the Vitali relation and its equivalence
classes have deep roots in descriptive set theory.
For instance Sierpinski \cite[c.~147]{sie}
and Luzin~\cite[Section~64]{lus:ea}
observed that the quotient set $\ves/\raz$ of all Vitali
classes has the property that (*) it cannot be mapped 
into $\ves$ by an injective Borel map.
On the other hand, as established in \cite{k91}, in
models of $\zf$ (without the axiom of choice) 
the Hartogs number of the set $\ves/\raz$  
(the least cardinal which cannot be injectively mapped
into $\ves/\raz$) 
can be greater than the continuum.
The relations $\Eo$ and $\vit$ play a key role in
modern studies of Borel equivalence relations,
being the least ones, 
in the sense of the Borel reducibility \cite{hkl}, 
among those satisfying (*).

\parf{Structure of the proof}
\las{ko8}

The proof of Theorem~\ref{Tun} is arganized as follows. 
Basic notions, related to Silver trees in the set  
$\bse$ of all finite dyadic strings, 
are introduced in sections \ref{tre}--\ref{nef}.  
Every set ${\jpi}$ of Silver trees $T$, closed under 
restriction to a given string $s\in T$, 
and under the natural action $s\app T$ by $s\in\bse$, 
is considered as \rit{a forcing by Silver trees}, 
a \rit{\sfo}, in brief. 
Every \sfo\ adjoins a \dd{\jpi}generic real $a\in\dn.$ 

Arguing in the constructible universe $\rL$, we  
define a forcing notion to prove Theorem~\ref{Tun} 
in Section~\ref{vop} 
in the form $\dP=\bigcup_{\al<\omi}\dP_\al$. 
The summands $\dP_\al$ are countable \sfo s defined 
by induction. 
Any \dd\dP generic extension of $\rL$ is a 
model for Theorem~\ref{Tun}. 
The inductive construction of $\dP_\al$ involves two 
key ideas.

The first idea, essentially by Jensen~\cite{jenmin}, 
is to make every level $\dP_\al$ of the construction  
\rit{generic} in some sense over the union of lower  
levels $\dP_\xi$, $\xi<\al$. 
This is based on a construction developed in  
sections \ref{gex} -- \ref{pres}, which includes   
the technique of fusion of Silver trees.
A special aspect of this construction, elaborated 
in \cite{kl22}, 
guarantees that $\dP$ 
is invariant under the group of transformations  
($2^{<\om}$ with the componentwise addition mod 2), 
which induces the equivalence relation $\Eo$. 
This invariance implies that $\dP$ 
(unless Jensen's original forcing) adjoins 
a \dd\Eo class of generic reals rather than a single 
such real as in \cite{jenmin}.
And overall, the successive genericity of the 
levels $\dP_\al$ implies that the three sets are 
equal in any \dd\dP generic extension of $\rL$: 

the \dd\Eo class $\eko a$ of the principal generic 
real $a\in\dn$, 

the intersection 
$\bigcap_{\al<\omi}\bigcup_{T\in\dP_\al}$, 

and the set of all \dd\dP generic reals over $\rL$. 

\noi
This equality, which leads to $\eko a$ being 
$\ip1\nn$, is established in 
sections \ref{pres+} -- \ref{bex2}.

The second idea goes back to old papers 
\cite{h74,ian78}.
In $\rL$, let $\vstf$ be the set of all countable 
sequences $\vjpi=\sis{\jpi_\xi}{\xi<\al}$ 
($\al<\omi$) 
compatible with the first genericity idea. 
Then a whole sequence $\sis{\jpi_\al}{\al<\omi}$ 
can be interpreted as a maximal branch in $\vstf$. 
It happens that if this branch is \rit{generic}, in 
some sense precisely defined in Section~\ref{vop} 
(\ref{ep2} of Theorem~\ref{ep}),  
with respect to all $\fs1{\nn-1}$ subsets of $\vstf$, 
then the ensuing forcing notion 
$\dP=\bigcup_{\al<\omi}\dP_\al$ 
inherits some basic forcing properties of the whole 
Silver forcing, up to a certain level 
of projective hierarchy. 
This includes, in particular, the invariance of the 
forcing relation with respect to some natural 
transformations of Silver trees, leading eventually 
to the proof of \ref{Tun3} of Theorem~\ref{Tun} 
in sections \ref{fs} -- \ref{zav}. 
But before that the actual construction of 
$\dP=\bigcup_{\al<\omi}\dP_\al$ 
involving both mentioned ideas is given 
in Section~\ref{vop}.

\parf{Silver trees} 
\las{tre}

Let $\bse$ be the set of all \rit{dyadic strings}   
\index{zz2om@$\bse$}%
\index{string!$\bse,$ dyadic strings}%
(finite sequences) 
\index{string}%
of numbers $0,1$ --- 
including  \rit{the empty string} $\La$.
\index{zzLa@$\La$, the empty string}%
\index{string!$\La$, the empty string}%
\index{zzLa@$\La$}%
If $t\in\bse$ and $i=0,1$, then  
$t\we i$ is the extension of $t$ by $i$ as the rightmost 
\index{zztiw@$t\we i$}%
term. 
If $s,t\in\bse$ then $s\sq t$ means that the string $t$ 
\index{zzstsq@$s\sq t$}%
\index{zzstsu@$s\su t$}%
extends $s$, while 
$s\su t$ means a proper extension. 
\index{string!extension, $\sq$, $\su$}%
The length of a string $s$ is denoted by $\lh s$,  
\index{zzlhs@$\lh s$, length}%
\index{string!length}%
\index{string!$2^n,$ strings of length $n$}%
and we let 
$2^n=\ens{s\in\bse}{\lh s=n}$ (strings of length $n$).%

As usual, $\dn$ is the set of 
all functions $f:\om\to2=\ans{0,1}$, 
--- \rit{the Cantor space}.
Any string $s\in\bse$ {\it acts\/} on $\dn$ so that 
$(s\app x)(k)=x(k)+s(k)\pmod 2$ whenever $k<\lh s$, and  
$(s\app x)(k)=x(k)$ otherwise.
\index{action!$s\app x$}%
\index{action!$s\app X$}%
\index{action!$s\app t$}%
If $X\sq\dn$ and $s\in\bse$ then let
$s\app X=\ens{s\app x}{x\in X}$.

Similarly, if $s\in 2^m,$ $t\in2^n,$ $m\le n$, then   
define a string $s\app t\in 2^n$ by 
$(s\app t)(k)=t(k)+s(k)\pmod 2$ whenever $k<m$, and  
$(s\app t)(k)=t(k)$ whenever $m\le k<n$.
If $m>n$, then let $s\app t=(s\res n)\app t$.
In both cases, $\lh{s\app t}=\lh t$. 
If $T\sq\bse$ then let $s\app T=\ens{s\app t}{t\in T}$.
\index{action!$s\app T$}%

\bdf
\lam{std}
A set  $T\sq\bse$ is a \rit{Silver tree}, 
\index{tree!Silver, $\pes$}%
\index{string!ukT@$u_k(T)$}%
in symbol $T\in\pes$, whenever there exists an infinite 
sequence of strings $u_k=u_k(T)\in\bse$ such that  
$T$ consists of all strings of the form  
$
s=u_0\we i_0\we u_1\we i_1\we u_2\we i_2 \we\dots\we u_{n}\we i_n
$, 
and their substrings (including $\La$), 
where $n<\om$ and $i_k=0,1$. 
 
In this case we let $\roo T=u_0$ 
({\it the stem\/} of $T$), 
\index{tree!stem, $\roo T$}%
\index{stem, $\roo T$}%
\index{zzstemT@$\roo T$}%
and define a closed set $[T]\sq\dn$  
\index{zzTII@$[T]$}%
\index{tree!branch}%
of all \rit{branches of $T$}, \ie\ all 
sequences  
$a=u_0\we i_0\we u_1\we i_1\we u_2\we i_2 \we\dots\in\dn,$ 
where $i_k=0,1$, $\kaz k$.
Let 
$$
\oin T{n}=\lh{u_0}+1+\lh{u_1}+1+\dots+\lh{u_{n-1}}+1+\lh{u_n}\,,
$$ 
in particular, $\oin T{0}=\lh{u_0}$, 
so that 
\index{zzsplnT@$\oin Tn$}%
$\oi T=\ens{\oin Tn}{n<\om}\sq\om$ is the set of all  
\index{zzsplT@$\oi T$}%
\rit{splitting levels} of $T\in\pes$; 
$\oi T$ is infinite.
If $u\in T\in\pes$, then 
define the \rit{restricted tree}  
\index{tree!restricted $T\ret u$}%
\index{zzTIu@$T\ret u$}%
$T\ret u=\ens{t\in T}{u\sq t\lor t\sq u}$; 
$T\ret u\in\pes$. 
\edf

\bpri
\lam{clop}
If $s\in\bse$ then the tree 
$\bad s=\ens{t\in\bse}{s\sq t\lor t\su s}$
\index{zzB:s:@$\bad s$}%
\index{tree!zzB:s:@$B[s]$}%
belongs to $\pes$, $\roo{\bad s}=u_0(\bad s)=s$, 
$\oin{\bad s}n=\lh s+n$, and if $k\ge1$ then   
$u_k(\bad s)=\La$.
In particular $\bad\La=\bse$. 
\epri

The corresponding sets  
$[T[s]]=\ens{x\in\bn}{s\su x}$ are clopen in $\dn$.

\ble
\lam{trans}
Let\/ $T\in\pes$.
Then\/ 
\ben
\renu
\itlb{trans01}%
if a set\/ $X\sq\dn$ is open and\/ $X\cap[T]\ne\pu$, 
then there is a string\/ $s\in T$ such that\/
$\req Ts\sq X\;;$ 

\itlb{trans1}%
if\/ $h\in\oi T$ and\/ $u,v\in2^h\cap T$ then\/
$\req Tv=(u\app v)\app\req Tu$ and\/ $(u\app v)\app T=T$. 
\een
\ele
\bpf
\ref{trans01}
If $a\in X\cap[T]$ then $\ans a=\bigcap_n\req T{a\res n}$, 
and hence 
$\req T{a\res n}\sq X$ for some $n$ by the compactness.
Let $s=a\res n$.
\epf

\parf{Splitting Silver trees}
\las{sst}

\rit{The simple splitting\/} of a tree $T\in\pes$ consists 
\index{splitting}%
of the subtrees
$$
\raw Ti=T\ret{\roo T\we i}=\ens{x\in[T]}{x(h)=i}\,,
\quad\text{where }\,h=\lh{\roo T},\;i=0,1.
$$
\index{splitting!zzT->i@$\raw T i$}%
\index{zzT->i@$\raw T i$}%
Then we have $\raw Ti\in\pes$, and  
$u_0(\raw Ti)=u_{0}(T)\we \ang{i}\we u_{1}(T)$, 
and then\linebreak[0]  
$u_k(\raw Ti)=u_{k+1}(T)$ whenever $k\ge1$, and 
$\oi{\raw Ti}= \oi T\bez\ans{\oin T0}$. 

The splitting can be iterated. 
Namely if
$s\in 2^n$ then define   
$$
\raw Ts=\raw{\raw{\raw{\raw T{s(0)}}{s(1)}}{s(2)}\dots}{s(n-1)}
\,.
$$
\index{splitting!zzT->s@$\raw T s$}%
\index{zzT->s@$\raw T s$}%
Separately let $\raw T\La=T$, for the empty string $\La$.

\ble
\lam{sadd}
Let\/ $T\in\pes$, $n<\om$, and\/ $h=\oin Tn$.
Then 
\ben
\renu
\itlb{sadd0}%
if\/ $s\in2^n$ then\/ $\raw Ts\in\pes$, 
$\lh{\roo{\raw Ts}}=h$, 
and there is a string\/
$u[s]\in 2^h\cap T$ such that\/ $\raw Ts=\req T{u[s]}\;;$ 

\itlb{sadd4}%
conversely if\/ $u\in 2^h\cap T$,
then there is a string\/ $s[u]\in2^n$ such that\/
$\req Tu=\raw T{s[u]}\;;$ 

\itlb{sadd4+}%
therefore, 
$\ens{\req Tu}{u\in T}=\ens{\raw Ts}{s\in\bse}\;.$ 
\een
\ele

\bpf
\ref{sadd0}
$u[s]=u_0(T)\we\ang{s(0)} 
\we \dots\we u_{n-1}(T)\we\ang{s(n-1)}\we u_{n}(T)$. 
\index{zzu.s@$u[s]$}%

\ref{sadd4}
Define $s=s[u]\in2^n$ by $s(k)=u(\oin Tk)$ 
for all $k<n$. 
\index{zzs.u@$s[u]$}%

\ref{sadd4+}
Let $u\in T$. 
Then $\oin T{n-1}<\lh u\le\oin Tn$ for some $n$. 
Now, by Definition~\ref{std}, 
there exists a (unique) string 
$v\in 2^h\cap T$, where $h=\oin Tn$, such that 
$\req Tu=\req Tv$. 
It remains to refer to \ref{sadd4}.
\epf

If $T,S\in\pes$ and $n\in\om$ then define $S\nq n T$
\imar{nq n}
(the tree $S$ \rit{\dd nrefines\/ $T$}), 
\index{zzSsubnT@$S\nq n T$}%
\index{zzsubn@$\nq n$}%
whenever  
$S\sq T$ and $\oin Tk=\oin Sk$ for all $k<n$.
In particular $S\nq 0 T$ is equivalent to just $S\sq T$.
By definition if $S\nq{n+1} T$ then $S\nq n T$
(and $S\sq T$).

\ble
\lam{tadd}
If\/ $T\in\pes$, $n<\om$, $s_0\in 2^n$, and\/ 
$U\in\pes$, $U\sq \raw T{s_0}$, then 
there is a unique tree\/ $T'\in\pes$ such that\/
$T'\nq n T$ and\/ $\raw{T'}{s_0}=U$. 
\ele

\bpf
Let $h=\oin Tn$. 
Pick a string $u_0=u[s_0]\in 2^h$ by 
Lemma~\ref{sadd}\ref{sadd0}.
Following Lemma~\ref{trans}\ref{trans1}, 
define $T'$ so that $T'\cap2^h=T\cap2^h$, and if
$u\in T\cap2^h$ then $\req {T'}u=(u\cdot u_0)\cdot U$.
In particular $\req{T'}{u_0}=U$.
\epf

\ble
\lam{fus}
Let\/
$\dots \nq 5 T_4\nq 4 T_3\nq 3 T_2\nq 2 T_1\nq 1 T_0$ 
be an infinite decreasing sequence of trees\/ $\pes$.
Then 
\ben
\renu
\itla{fus1}
$T=\bigcap_nT_n\in\pes\;;$ 
\itla{fus2}
if\/ $n<\om$ and\/ $s\in2^{n+1}$ then\/
$\raw Ts=T\cap\raw{T_{n}}s=\bigcap_{m\ge n}\raw{T_m}s$.
\een
\ele
\bpf
Note that  
$\oi T=\ens{\oin{T_n}{n}}{n<\om}$; both claims  
easily follow.
\epf

The next application of this lemma 
and Lemma~\ref{tadd} is well known:

\bcor
\lam{ssq}
Let\/ $T\in\pes$. 
If a set\/ $X\sq[T]$ has the Baire property inside\/ $[T]$, 
but is not meager inside\/ $[T]$, 
then there is a tree\/ $S\in\pes$ such that\/ 
$[S]\sq X$. 
If\/ $f:[T]\to\dn$ is a continuous map then there exists  
a tree\/ $S\in\pes$ such that\/ $S\sq T$ and\/ $f\res[S]$ 
is a bijection or a constant.\qed
\ecor

\parf{\sfo s}
\las{sitf}

\bdf
\lam{fds}
\index{forcing!by Silver trees}%
\index{forcing!\sfo}%
Any set $\jpi\sq\pes$ satisfying  
\ben
\Aenu
\itlb{sitf2}%
if $u\in T\in{\jpi}$ then $\req Tu\in {\jpi}$, and

\itlb{sitf3}%
if $T\in{\jpi}$ and $\sg\in\bse$ then $\sg\app T\in {\jpi}$, 
\een
is called {\it a forcing by Silver trees}, 
\rit{\sfo} in brief. 
\edf

\bre
\lam{xg}
Any \sfo\ ${\jpi}$ can be considered as a forcing notion  
ordered so that if $T\sq T'$, then $T$ is a stronger
\usl. 
The forcing ${\jpi}$ adjoins a real $x\in\dn$. 
More exactly if a set $G\sq{\jpi}$ is 
\dd{\jpi}generic over a given model or set universe 
$M$ (and ${\jpi}\in M$ is assumed)  
then the intersection $\bigcap_{T\in G}[T]$ contains a  
unique real $a=\xx G \in \dn$, and this real satisfies  
\index{zzaG@$\xx G$, generic real}%
\index{generic real!$\xx G$}%
$M[G]=M[\xx G]$ and $G=\ens{T\in{\jpi}}{\xx G\in[T]}$. 
Reals $\xx G$ of this kind are called  
\dd{\jpi}\rit{generic}. 
\ere

\bpri
\lam{xcoh}
The following sets are \sfo s: 
the set  
$\dpo=\ens{T[s]}{s\in\bse}$ 
of all trees in \ref{clop} ---  
\rit{the Cohen forcing}, 
\index{forcing!Cohen forcing, $\dpo$}%
\index{forcing!Silver forcing, $\pes$}%
\index{zzPcoh@$\dpo$}
and the set $\pes$ of all Silver trees 
--- the \rit{Silver forcing} itself.
\epri

\ble
\lam{komr}
If\/ $\pu\ne\rQ\sq\pes$ then the following set
is a \sfo$:$ 
$$
{\jpi}
=\ens{\sg\app(\req Tu)}{u\in T\in\rQ\land\sg\in\bse}
=\ens{\sg\app(\raw Ts)}{T\in\rQ\land s,\sg\in\bse}.
$$
\ele
\bpf
To prove \ref{fds}\ref{sitf2}, 
let $T\in\rQ$ and $v\in S=\sg\app(\req Tu)$. 
Then $w=\sg\app v\in \req Tu$ and $v=\sg\app w$. 
It follows that $\req Sv=\sg\app(\req{\req Tu}w)$, 
where $\req{\req Tu}w=\req Tu$ or $=\req Tw$, whenever 
accordingly $w\sq u$ or $u\su w$.
The second equality of the lemma   
follows from Lemma~\ref{sadd}\ref{sadd4+}.
\epf

\bdf
[collages] 
If ${\jpi}\sq\pes$, $T\in\pes$, $n<\om$, and all split  
trees $\raw Ts$, $s\in 2^n$, belong to ${\jpi}$ then   
$T$ is an \rit{\dd ncollage over\/ ${\jpi}$}.
\index{collage}%
The set of all \dd ncollages over ${\jpi}$ is  
$\sct{\jpi} n$.  
\index{zzColgnP@$\sct{\jpi} n$}%
\imar{sct Pn}%
Then 
${\jpi}=\sct{\jpi}0\sq\sct{\jpi} n\sq \sct{\jpi}{n+1}$. 
\edf

\ble
\lam{stf}
Let\/ ${\jpi}\sq\pes$ be a \sfo, and\/ $n<\om$.
Then$:$
\ben
\renu
\itlb{nq2}%
if\/ $T\in {\jpi}$ and\/ $s\in\bse$ then\/ 
$\raw Ts\in{\jpi}\;;$ 


\itlb{stf1}%
if\/ $T\in\pes$ and\/ $s_0\in2^n$ 
then\/ $\raw T{s_0}\in{\jpi}$ is equivalent to\/ 
$T\in\sct{\jpi} n$$;$ 

\itlb{stf2}%
if\/ $U\in\sct{\jpi} n$, $s_0\in2^n$, $S\in{\jpi}$, and\/
$S\sq \raw U{s_0}$ then there is a tree\/ 
$V\in\sct{\jpi} n$ such that\/ $V\nq n U$ and\/
$\raw{V}{s_0}=S\;;$ 

\itlb{stf2c}%
if\/ $U\in\sct{\jpi} n$ and\/ 
$D\sq{\jpi}$ is open dense in\/ ${\jpi}$ then there is
a tree\/ $V\in\sct{\jpi} n$ such that\/ $V\nq n U$
and\/ $\raw{V}{s}\in D$
for all\/ $s\in 2^n$.
\een
\ele

A set\/ $D\sq{\jpi}$ is \rit{dense} in ${\jpi}$  
if for any $S\in{\jpi}$ there is a tree $T\in D$, 
$T\sq S$, and \rit{open dense}, if in addition   
$S\in D$ holds whenever $S\in{\jpi}$, $T\in D$, $S\sq T$.

\bpf
To prove \ref{nq2} make use of \ref{fds}\ref{sitf2} 
and Lemma~\ref{sadd}\ref{sadd0}.

\ref{stf1}
If $T\in\sct{\jpi} n$ then by definition $\raw T{s_0}\in{\jpi}$.
To prove the converse let $h=\oin Tn$ and let a string 
$u[s]\in 2^h\cap T$ satisfy 
$\raw Ts=\req T{u[s]}$ for all $s\in2^n$
by lemma~\ref{sadd}\ref{sadd0}.
If $\raw T{s_0}\in{\jpi}$ and $s\in 2^n$ then 
$\raw Ts=\req T{u[s]}=(u[s]\app u[s_0])\app \req T{u[s_0]}
=(u[s]\app u[s_0])\app \raw T{s_0}$
by Lemma~\ref{trans}\ref{trans1}.
Therefore $\raw Ts\in{\jpi}$ 
by \ref{fds}\ref{sitf3}.
And finally $T\in \sct{\jpi} n$.

\ref{stf2}
By Lemma~\ref{tadd}, there is a tree $V\in\pes$,
satisfying $V\nq n U$ and $\raw{V}{s_0}=S$.
However $V\in\sct{\jpi} n$ by \ref{stf1}.

To prove \ref{stf2c}
apply \ref{stf2} $2^n$ times (for all $s\in2^n$). 
\epf

\parf{Continuous maps} 
\las{nef}

Let ${\jpi}$ be a fixed \sfo\ in this section. 

{\ubf Regularity.} 
We study the behaviour of continuous maps on sets  
of the form $[T]$, $T\in  {\jpi}$. 
The next definition highlights the case when a given 
continuous $f:\dn\to\dn$ is forced to be not equal 
to a map of the form $x\mto\sg\app x$, $\sg\in\bse.$ 

\bdf
\lam{ref}
Let $T\in {\jpi}$.
A continuous $f:\dn\to\dn$ is
\index{regular}%
\index{map!regular}%
\index{function!regular}%
\rit{regular on\/ $T$ inside\/ ${\jpi}$}, 
if there is no tree $T'\in {\jpi}$ and string $\sg\in\bse$ 
such that $T'\sq T$ and $\sg\app f(x)=x$
(equivalently, $f(x)=\sg\app x$) 
for all $x\in[T']$.
\edf

\ble
\lam{rc}
Let\/ $S,T\in {\jpi}$, $f:\dn\to\dn$
be continuous, 
and\/ $\sg\in\bse$.
Then$:$
\ben
\renu
\itla{rc1}
there are trees\/ $S',T'\in {\jpi}$ such that\/ 
$S'\sq S$, $T'\sq T$,   
$[T']\cap(\sg\app\ima f[S'])=\pu\;;$ 

\itla{rc2}
if\/ $\tau\in\bse$, $T=\tau\app S$, and\/ $f$ 
is regular on\/ $S$ inside\/ ${\jpi}$, 
then there exist trees\/ $S',T'\in {\jpi}$ such that\/ 
$S'\sq S$, $T'\sq T$, 
$T'=\tau\app S'$, and\/ 
$[T']\cap(\sg\app\ima f[S'])=\pu\;.$ 
\een
\ele
\bpf
\ref{rc1}
Let $x_0\in[S]$ and let $y_0\in[T]$ be different  
from $\sg\app f(x_0)$. 
As $f$ is continuius, there is $m\ge\lh\sg$, 
such that 
$\sg\app(f(x)\res m)\ne y_0\res m$ whenever  
$x\in[S]$ and $x\res m=x_0\res m$. 
By \ref{fds}\ref{sitf2}, 
there are trees\/ $S',T'\in {\jpi}$ such that  
$S'\sq S$, $T'\sq T$,  
$x\res m=x_0\res m$ for all $x\in[S']$, and
$y\res m=y_0\res m$ for all $y\in[T']$.\vom

\ref{rc2}
Assume that $\lh\sg=\lh\tau$ 
(otherwise the shorter string can be extended by zeros). 
The set $X=\ens{x\in[S]}{\sg\app f(x)\ne \tau\app x}$ 
is open in $[S]$ and non-empty, by the regularity. 
Let $x_0\in X$. 
There is a number $m\ge\lh\sg=\lh\tau$ satisfying  
$\tau\app(x_0\res m)\ne \sg\app(f(x_0)\res m)$. 
By \ref{fds}\ref{sitf2} there is  
a tree $S'\in {\jpi}$ satisfying $x_0\in[S']$ and  
$x\res m=x_0\res m$, $f(x)\res m=f(x_0)\res m$ 
for all $x\in[S']$.
Put $T'=\tau\app S'$. 
\epf

{\ubf Coding continuous maps.}  
If $f:\dn\to\dn$ is continuous, $k<\om$ and $i=0,1$
then the set $\bfa fki=\ens{x\in\dn}{f(x)(k)=i}$,
and  
$\bfa fk0\cap\bfa fk1=\pu$,  
$\bfa fk0\cup\bfa fk1=\dn$.
Let $\kod fki$ 
be the set of all \dd\sq minimal strings $s\in\bse$  
such that $\bad s\sq \bfa fki$.
Then $\kod fk0,\kod fk1\sq\bse$ are finite  
disjoint irreducible antichains\snos
{An antichain $A\sq\bse$ 
is \rit{irreducible}, if it does not contain a  
pair of the form $s\we 0,s\we 1$.} 
in $\bse$, and 
$\bfa fki=\bigcup_{s\in\kod fki}\bad s$, 
while the union $\kod fk0\cup\kod fk1$ is a  
maximal antichain in $\bse.$ 
Let $\koe f=\sis{\kod fki}{k<\om,i=0,1}$ 
--- the \rit{code} of $f$.
\index{function!code, $\koe f$}%
\index{code, $\koe f$}%
\index{zzcode@$\koe f$}%

The other way around, if
$\loe=\sis{\lod ki}{k<\om,i=0,1}$ is a family  
of finite irreducible antichains $\lod ki\sq\bse$,
and if $k<\om$ then  
$\lod k0\cap\lod k1=\pu$ while  
$\lod k0\cup\lod k1$ is a maximal antichain in $\bse$, 
then $\loe$ is called \rit{a code of continuous function}. 
In this case, the related continuous function  
$f:\dn\to\dn$ is denoted by $f_\loe$, that is, 
$f_\loe(x)(k)=i$ whenever  
$x\in\bigcup_{s\in\lod ki}\bad s$.

Let $\ccf$ denote the set of all
codes of continuous functions.

\parf{Generic extensions of \sfo s} 
\las{gex}

The forcing notion to prove Theorem~\ref{Tun}
will be defined in the form of an \dd\omi union 
of its countable parts --- \rit{levels}.
The next definition presents requirements which
will govern the interactions between the levels.

\bdf
\lam{dfex}
Let $\cM$ be any set and $\jpi$ be a \sfo.  
We say that another \sfo\ $\jqo$ 
is an \dd\cM\rit{extension} of $\jpi$, 
in symbol $\jpi\prol\cM\jqo$, if the following holds:%
\index{forcing!\dd\cM extension, $\prol\cM$}%
\index{Mextension@\dd\cM extension, $\prol\cM$}%
\index{zzsqsubM@$\prol\cM$}%
\index{zzsqsub@$\prok$}%
\imar{prol cM}%
\ben
\Aenu
\itlb{dfex1}%
$\jqo$ is dense in $\jqo\cup\jpi$; 

\vyk{
\itlb{dfex2x}%
if $U\in\jqo$ and $T\in\jpi$ then $[U]\cap[T]$ 
is clopen in $[U]$;
}

\itlb{dfex3}%
if a set $D\in\cM$, $D\sq\jpi$ is pre-dense in $\jpi$ 
and $U\in\jqo$ then $U\sqf\bigcup D$, that is,
\index{zzsqfin@$\sqf$}%
there is a finite set $D'\sq D$ such that 
$U\sq\bigcup D'$;

\vyk{
\itlb{dfex4}%
if $T,T'\in\jpi$   
are incompatible in $\jpi$ then $T,T'$  
are incompatible in $\jqo\cup\jpi$, too, and moreover  
if $U\in\jqo$ then $[T]\cap[T']\cap[U]=\pu\;;$
}

\itlb{dfex6}%
if $T_0\in\jpi$,    
$f:\dn\to\dn$ is a continuous map with a code in  
$\cM$, regular on $T_0$ inside ${\jpi}$, and $U,V\in\jqo$, 
$U\sq T_0$, 
then $[V]\cap(\ima f[U])=\pu$.\qed 
\een
\eDf

Speaking of \ref{dfex6}, it claims essentially that  
if the regularity holds then  
$T_0$ forces, in the sense of
$\jqo\cup\jpi$, that $f(\xx G)$ does not belong
to $\bigcup_{V\in\jqo}[V]$.

If $\cM=\pu$ then we write $\jpi\prok\jqo$ 
\index{forcing!extension, $\prok$}%
\imar{prok}%
instead of $\jpi\prol\pu\jqo$; 
in this case \ref{dfex3} and \ref{dfex6}
are trivial.
Generally, we'll consider, in the role of $\cM$, 
transitive models of the theory $\zfc'$ 
\index{zzzfcp@$\zfc'$}%
\index{theory!zfcp@$\zfc'$}%
which includes all $\ZFC$ axioms except for the Power   
Set axiom, but an axiom is adjoined, which claims  
the existence of $\pws\om$. 
(Then the existence of  $\omi$ and sets like $\dn$ 
and $\pes$ 
easily follows.)

\ble
\lam{tec}
Let\/ $\cM\models\zfc'$ be a transitive model,    
$\jpi\in \cM$ and\/ $\jqo$ be\/ 
\sfo s satisfying\/ $\jpi\prol\cM\jqo$, and\/ 
$U\in\jqo$.  
Then 
\ben
\renu
\itlb{tec1}%
if\/ $T\in\jpi$ then\/ $[U]\cap[T]$ 
is clopen in\/ $[U]\,;$ 

\itlb{tec2}%
if\/ $D\in\cM$, $D\sq\jpi$ is pre-dense in\/ $\jpi$  
then\/ $D$ remains pre-dense in $\jqo\cup\jpi\,;$

\itlb{tec4}%
if\/ $T,T'\in\jpi$   
are incompatible in\/ $\jpi$ then\/ $T,T'$  
are incompatible in\/ $\jqo\cup\jpi$, too, and moreover  
if\/ $U\in\jqo$ then\/ $[T]\cap[T']\cap[U]=\pu\;;$ 

\itlb{tec3}%
if\/ $R\sq U$ is a Silver tree 
then the set\/ 
$\bR=\jqo\cup\ens{\sg\app (\raw Rt)}{t,\sg\in\bse}$ 
is a \sfo, and still\/ 
$\jpi\prol\cM\bR$. 
\een
\ele
\bpf
\ref{tec1}
The set 
$D(T)=\ens{S\in\jpi}{S\sq T\lor [S]\cap[T]=\pu}$
belongs to $\cM$ and is open dense in $\jpi$ by 
\ref{fds}\ref{sitf2}.
Then by \ref{dfex}\ref{dfex3} there is a finite 
set $D'\sq D(T)$ such that $U\sq\bigcup D'$. 
But then we have 
$[U]\bez [T]=\bigcup_{S\in D''}([S]\cap[W])$, 
where $D''=\ens{S\in D'}{[S]\cap[T]=\pu}$. 
Thus $[U]\bez[T]$ is a closed set.

\ref{tec2}
Let $U\in\jqo$. 
Then by \ref{dfex}\ref{dfex3} there is a finite 
set $D'\sq D$ such that $U\sq\bigcup D'$.
There is a tree $T\in D'$ such that $[T]\cap [U]$ 
has a non-empty open interior in $[U]$.
By Lemma~\ref{trans}\ref{trans01}, 
there is $s\in U$ such that   
$[U']\sq [T]\cap[U]$, where $U'=\req Us$, 
thus $U'\sq U\cap T$. 
Finally $U'\in\jqo$, as $\jqo$ is a \sfo.

\ref{tec4}
By the incompatibility, 
if $S\in\jpi$ then $S\not\sq T$ or $S\not\sq T'$, 
and hence there is a tree $S'\in\jpi$, 
$S'\sq S$, satisfying $[S']\cap[T]\cap [T']=\pu$. 
Therefore the set 
$D=\ens{S\in\jpi}{[S]\cap[T]\cap[T']=\pu}$ 
is dense in $\jpi$ and belongs to $\cM$. 
Now if $U\in\jqo$ then $U\sqf\bigcup D$ by 
\ref{dfex}\ref{dfex3}, thus immediately 
$[U]\cap[T]\cap[T']=\pu$.

\ref{tec3}
$\bR$ is a \sfo\ by Lemma~\ref{komr}, so we have 
to prove $\jpi\prol\cM\bR$.  
We skip a routine verification of 
\ref{dfex}\ref{dfex1},\ref{dfex3} 
and focus on \ref{dfex6}. 
Let $T_0$, $f$ be as in \ref{dfex6}. 

Consider trees  
$R'=\sg\app (\raw Rt)\in\dR$, $R'\sq T_0$, and 
$V\in\jqo$; we have to prove that 
$[V]\cap(\ima f[R'])=\pu$.
The problem is that while 
$R'\sq\sg\app R\sq U'=\sg\app U$, 
it may not be true that $U'\sq T_0$. 
But, now we claim that 
\rit{there is a finite set of trees\/ 
$\ans{U_1,\dots,U_n}\sq\jqo$ such that\/ 
$R'\sq U_1\cup\ldots\cup U_n\sq T_0$.} 
If this is established then  
$[V]\cap(\ima f[U_i])=\pu$ for all $i$ since 
$\jpi\prol\cM\jqo$, and we are done.  

To prove the claim consider the dense set 
$D(T_0)\in\cM$ as in the proof of \ref{tec1} above.
Then $U'\sqf\bigcup D$, hence there is a finite 
set $D'\sq D(T_0)$ satisfying $U'\sq\bigcup D'$. 
As by construction $R'\sq U'\cap T_0$, 
we conclude that $R'\sq\bigcup D''$, 
where $D''=\ens{S\in D'}{S\sq T_0}$. 
Now if $S\in D''$ then $[S]\cap [U']$ is clopen 
in $[U']$, 
therefore $[S]\cap [U']=\bigcup_{v\in W(S)}\req {U'}v$, 
where $W(S)\sq U'$ is a suitable finite set. 
We put 
$\ans{U_1,\dots,U_n}=
\bigcup_{S\in D''}\ens{\req {U'}v}{v\in W(S)}$.
\epf 

\bte
\lam{tex}
Let\/ $\cM\models\zfc'$ be a countable transitive model    
and\/ $\jpi\in \cM$ be a \sfo.  
Then there exists a countable\/ \sfo\ $\jqo$
satisfying\/ $\jpi\prol\cM\jqo$.
\ete

The proof of this theorem is presented below. 
Section~\ref{jex} contains the definition of 
the \sfo\ $\jqo$, while the proof of its 
properties follows in Section \ref{pres}.

\parf{Construction of extending \sfo} 
\las{jex}

The next definition formalizes a construction of   
countably many Silver trees by means of  
Lemma~\ref{fus}.   

\bdf
\lam{muss}
A \rit{\mus} is any indexed set 
\index{system}%
\index{system!$\abs\vpi$}%
$\vpi=
\sis{\ang{\yh\vpi m,\xt\vpi m}}{m\in\abs\vpi}$, 
\index{system!$\zh\vpi\xi m$}%
\index{system!$\zt\vpi\xi m$}%
\index{zznfxm@$\zh\vpi\xi m$}%
\index{zztfxm@$\zt\vpi\xi m$}%
where $\abs\vpi\sq \om$ is finite, 
and if $m\in\abs\vpi$ then
$\yh\vpi m\in\om$, 
\index{zzfII@$\abs\vpi$}%
\index{support!$\abs\vpi$}%
$\xt\vpi m=
\ang{\xT\vpi m0,\xT\vpi m1,\dots,
\xT\vpi m{\yh\vpi m}}$, 
each $\xT\vpi mn$ is a Silver tree, 
\index{system!$\zT\vpi\xi mn$}%
\index{zzTfxmn@$\zT\vpi\xi mn$}%
\index{tree!zzTfxmn@$\zT\vpi\xi mn$}%
and $\xT\vpi  m{n+1}\nq {n+1} \xT\vpi m{n}$ 
whenever $n<\yh\vpi m$. 

In this case, if $n\le \yh\vpi m$ and $s\in2^n$ 
then let $\xT\vpi ms=\raw{\xT\vpi mn}{s}$. 
\index{system!$\zT\vpi\xi ms$}%
\index{zzTfxms@$\zT\vpi\xi ms$}%
\index{tree!zzTfxms@$\zT\vpi\xi ms$}%

If $\jpi$ is a \sfo\ then $\ms\jpi$ 
\index{system!$\ms\jpi$}%
\index{zzMSP@$\ms\jpi$}%
is the set of all \mus\ $\vpi$ such that  
$\xT\vpi mn\in\sct{\jpi}n$ for all  
$m\in\abs\vpi$ and $n\le \yh\vpi m$.
Then every tree $\xT\vpi ms$ 
belongs to $\jpi$.
\edf

A \mus\  
$\vpi$ \rit{extends} a \mus\ $\psi$, 
\index{system!extension, $\psi\cle\vpi$}%
in symbol $\psi\cle\vpi$, 
\imar{cle}%
if $\abs \psi\sq\abs\vpi$, and  
for any index $m\in\abs\psi$, first,  
$\yh\vpi m\ge \yh\psi m$, and second,   
$\xt\vpi m$ extends $\xt\psi m$ in the sense that    
simply $\xT\vpi m{n}=\xT\psi m{n}$ for all   
$n\le\yh\psi m$.

\ble
\lam{nwm}
Let\/ $\jpi$ be a \sfo\ and\/ 
$\vpi\in\ms\jpi$. 

If\/ $m\in\abs\vpi$ and\/ $n=\yh\vpi m$ then the
extension\/ $\vpi'$ of the \mus\/ $\vpi$ by\/
$\yh{\vpi'} m=n+1$ and\/
$\xT{\vpi'} m{n+1}=\xT{\vpi} m{n}$ belongs to\/
$\ms\jpi$ and\/ $\vpi\cle\vpi'$.  

If\/ $m\nin\abs\vpi$ and\/ $T\in\jpi$ then the
extension\/ $\vpi'$ of the \mus\/ $\vpi$ by\/
$\abs{\vpi'}=\abs\vpi\cup\ans{m}$, 
$\yh{\vpi'} m=0$ and\/
$\xT{\vpi'} m{0}=T$,  
belongs to\/ $\ms\jpi$ and\/ $\vpi\cle\vpi'$.\qed 
\ele

Now, according to the formulation of Theorem~\ref{tex}, 
we assume that 
$\cM\models\zfc'$ is a countable transitive model 
and $\jpi\in\cM$ is a (countable) \sfo.

\bdf
\lam{dPhi}
(i) 
We fix a \dd\cle increasing sequence   
$\dphi=\sis{\vpi(j)}{j<\om}$ 
of \mus s   
$\vpi(j) \in\ms\jpi$, 
\rit{generic over\/ $\cM$} in the sense that it 
intersects every set 
$D\in\cM\yd D\sq\ms\jpi$ dense in $\ms\jpi$.
The density here means that for any \mus\ $\psi\in\ms\jpi$ 
there is a \mus\ $\vpi\in D$ such that $\psi\sq\vpi$.  
The goal of this definition is to define another 
\sfo\ $\jqo$ from $\dphi$, see item (iv) below. \vom

(ii) 
If $m,n<\om$ then the set  
$ 
D_{m n}=\ens{\vpi\in\ms\jpi}
{\yh\vpi m\ge n}  
$ 
is dense by Lemma~\ref{nwm} and belongs to $\cM$, 
hence it intersects $\dphi$.
Therefore if $m<\om$ then there exists an infinite 
sequence  
$$
\dots \nq 5 \tx{m}{4}\nq 4 \tx{m}{3} 
\nq 3 \tx{m}{2}\nq 2 \tx{m}{1}\nq 1 \tx{m}{0}
$$
of trees $\tx{m}{n} \in\sct{\jpi}n$ such that 
\index{zzTFxmn@$\tx{\xi m}{n}$}%
\index{tree!zzTFxmn@$\tx{\xi m}{n}$}%
for each $j$, if 
$m\in\abs{\vpi(j)}$ and $n\le\yh{\vpi(j)}m$, 
then $\ntt{\vpi(j)}{m}n = \tx{m}{n}$.
If $n<\om$ and $s\in2^n$ then let 
$\tx{m}s=\raw{\tx{m}n}s$; 
\index{zzTFxms@$\tx{\xi m}{s}$}%
\index{tree!zzTFxms@$\tx{\xi m}{s}$}%
then $\tx{m}s\in \jpi$ because  
$\tx{m}n\in\sct{\jpi}n$.\vom

(iii)
In this case, by Lemma~\ref{fus}, every set 
$$
\TS
\ufi{m}
=\bigcap_{n}\tx{m}n=
\bigcap_{n}\bigcup_{s\in2^n}\tx{m}s
$$
\index{zzUFxm@$\ufi{\xi m}$}%
\index{tree!zzUFxm@$\ufi{\xi m}$}%
belongs to $\pes$, 
the same is true for all subtrees $\qfs{}ms$, 
and we have 
$$
\TS
\qfs{}ms=\ufi{m}\cap \tx{m}s=
\bigcap_{n\ge \lh s}{\raw{\tx{m}n}s} 
\,,
$$
by Lemma~\ref{fus}, 
and it is clear that $\ufi{m}=\qfs{}m\La$. 
In addition if strings $s, t\in\bse$ 
satisfy\/ $s\sq t$ then  
$\ty{}ms\sq\ty{}mt$ and $\qfs{}ms\sq\qfs{}mt$, 
but if $s,t$ are \dd\sq incomparable then  
$[\qfs{}ms]\cap[\qfs{}mt]=
[\ty{}m{s}] \cap[\ty{}mt]=\pu$.\vom

(iv) 
Finally the set    
$\jqo=\ens{\sg\app\qfs{}ms}{m<\om\land \sg,s\in\bse}$ 
is a countable \sfo\ by Lemma~\ref{komr}. 
\edf

\parf{Validation of the extension property} 
\las{pres}

Here we prove that $\jpi\prol\cM\jqo$ 
in the context of Definition~\ref{dPhi}. 
We check all requirements of Definition \ref{dfex}.\vom

{\bf \ref{dfex}\ref{dfex1}.}
If $T\in\jpi$ then  
the set $\Da(T)$ of all 
\mus s $\vpi\in\ms\jpi\cap\cM$ 
such that $\zT\vpi{}m0=T$ for some $m$, 
belongs to $\cM$ and is dense in $\ms\jpi$ by 
Lemma~\ref{nwm}. 
Therefore $\vpi(J)\in \Da(T)$ for some $J$, 
by the choice of $\dphi$. 
Then $\ty{}m0=T$ for some $ m<\om$. 
But $\qfs{}m\La=\qfd{}m \sq \ty{}m0$ and 
$\qfd{}m\in\jqo{}$.\vom

{\bf \ref{dfex}\ref{dfex3}.}
Let $U=\sg\app\qfs{}ms\in\jqo$, 
where $m<\om$ and $s,\sg\in\bse.$  
We have to check $U\sqf\bigcup D$. 
We can assume, as above, that $\sg=\La$, 
that is, $U=\qfs{}ms$. 
(Otherwise substitute $\sg\app D$ for $D$.) 
We can also assume that $s=\La$, 
that is, $U=\qfd{}m$, as  
$\qfs{}ms\sq\qfd{}m$.
Thus let $U=\qfd{}m$. 

The set $\Da\in\cM$ of all \mus s  
$\vpi\in \ms\jpi$ such that  
$m\in\abs\vpi$ and for every string $t\in 2^{n}$, 
where $n=\zh\vpi{}m$, there is a tree $S_t\in D$ with   
$\zT\vpi{}m t\sq S_t$,   
is dense in $\ms\jpi$ by lemma \ref{stf}\ref{stf2c} 
due to the pre-density of $D$ itself. 
Therefore there is an index $j$ such that 
$\vpi(j)\in\Da$.
Let this be witnessed by trees  
$S_t\in D\yt t\in 2^{n},$ 
where $n=\zh{\vpi(j)}{}m$, 
so that $\zT{\vpi(j)}{}m t\sq S_t$, $\kaz t$, 
and hence $\zT{\vpi(j)}{}m n\sqf D$. 
Then 
$ 
\textstyle
U\,=\,\qfd{}m
\,\sq\,
\tfs{}mn
\,=\,
\zT{\vpi(j)}{}m n 
\,\sqf\,
\bigcup D\,.
$  
\vom

\vyk{
{\bf \ref{dfex}\ref{dfex2x}.} 
Let $T\in\jpi$ and $U=\qfs{}ms\in\jqo$.
(If $U=\sg\app\qfs{}ms$, $\sg\in\bse,$ then  
replace $T$ by $\sg\app T$.)
The set $D$ of all trees $W\in\jpi$ 
such that $W\sq T$ or $[W]\cap[T]=\pu$, 
is clopen in $\jpi$ by \ref{fds}\ref{sitf2}. 
Therefore we have $U\sqf\bigcup D$ 
by \ref{dfex}\ref{dfex3}, thus $U\sqf\bigcup D'$, 
where $D'\sq D$ is finite. 
But then $[U]\bez [T]=\bigcup_{W\in D''}([U]\cap[W])$, 
where $D''=\ens{W\in D'}{[W]\cap[T]=\pu}$. 
Thus $[U]\bez[T]$ is a closed set.\vom
}

\vyk{
{\bf \ref{dfex}\ref{dfex4}.}
Assume that $T,T'\in\jpi$ are incompatible in $\jpi$. 
Then if $S\in\jpi$ then $S\not\sq T$ or $S\not\sq T'$, 
and hence there is a tree $S'\in\jpi$, 
$S'\sq S$, satisfying $[S']\cap[T]\cap [T']=\pu$. 
Therefore the set $D$ of all trees $S\in\jpi$ such  
that $[S]\cap[T]\cap [T']=\pu$, is dense in $\jpi$. 
Now apply Lemma~\ref{tec}\ref{tec2}.\vom 
} 

{\bf \ref{dfex}\ref{dfex6}.}
Let $T_0\in\jpi$, $f$, and $U,V\in\jqo$ be as in  
\ref{dfex}\ref{dfex6}. 
Then $U=\tau\app\qfs{}K{s_0}$, where  
$s_0,\tau\in\bse$ and $K<\om$. 
We can wlog assume that $\tau=\La$, 
that is, $U=\qfs{}K{s_0}$, since the general case is 
reducible to this case by the substitution 
of $\tau\app T_0$ for $T_0$ and  
the function $f'(x)=f(\tau\app x)$ for $f$. 
Thus let $U=\qfs{}K{s_0}$. 

Similarly, generally speaking 
$V=\rho\app\qfs{}L{t_0}$, where $t_0,\rho\in\bse$ 
and $L<\om$. 
But this is reducible to the case $\rho=\La$  
by the substitution of  
$f'(x)=\rho\app f(x)$ for $f$.
Thus let $V=\qfs{}L{t_0}$. 
Finally, since $\qfs{}L{t_0}\sq \qfd{}L$, we can  
assume that even $V=\qfd{}L$.
Now consider the set $\Da\in\cM$ of all \mus s   
$\vpi\in \ms\jpi$ such that there is 
a number $m<\om$ satisfying the following:
\ben
\Renu
\itlb{ww42}\msur
$K,L\in\abs\vpi$, $\yh\vpi K=\yh\vpi L=m$, 
and $\lh{s_0}\le m$; 

\itlb{ww41}%
if $s\in2^m$ then $\xT\vpi Ks \sq T_0$ or  
$[\xT\vpi Ks]\cap [T_0]=\pu$; 

\itlb{ww43}%
if $s\in2^m$ and $\xT\vpi Ks \sq T_0$ then 
$\xT\vpi Lm\cap (\ima f[\xT\vpi Ks])=\pu$. 
\een

\ble
\lam{Ka}
The set\/ $\Da$ is dense in\/ $\ms\jpi$. 
\ele
\bpf
Let $\psi\in\ms\jpi$; 
we have to define a \mus\ $\vpi\in\Da$ satisfying   
$\psi\cle\vpi$.
By lemma~\ref{nwm}, we can assume that $K,L\in\abs\psi$ 
and  
$\yh\psi K=\yh\psi L=m-1$ for some $m\ge\lh{s_0}$.
Now we define a first preliminary version of  
$\vpi$, by $\yh\vpi K=\yh\psi L=m$ and 
$\xT\vpi Km=\xT\psi K{m-1}$, $\xT\vpi Lm=\xT\psi L{m-1}$, 
and keeping the other elements of $\vpi$ equal to  
those of the system $\psi$, so that $\psi\cle\vpi$. 
The trees $S=\xT\vpi Km$ and $T=\xT\vpi Lm$ 
belong to $\sct \jpi m$. 

The set $D(T_0)$ of all trees $W\in\jpi$, 
such that $W\sq T_0$ or $[W]\cap[T_0]=\pu$, 
is open dense in $\jpi$ by \ref{fds}\ref{sitf2}. 
Therefore by Lemma~\ref{stf}\ref{stf2c}    
there is a tree $S'\in\sct \jpi m$ satisfying  
$S'\nq m S$ and if $s\in2^m$ then $\raw {S'}s\in D(T_0)$.
We modify $\vpi$ by putting $\xT\vpi Km=S'$
instead of $S$. 
It is clear that the new system still satisfies  
$\psi\cle\vpi$, and in addition \ref{ww41} 
holds by construction.

Further modification of $\vpi$ towards \ref{ww43} 
depends on whether $K=L$.

{\ubf Case 1:} $K\ne L$. 
If $s,t\in2^m$ then by Lemma~\ref{rc}\ref{rc1}  
there exist trees $C,D\in\jpi$ such that 
$C\sq\raw Ss$, $D\sq\raw Tt$, and  
$[D]\cap (\ima f{[C]})=\pu$. 
Applying Lemma~\ref{stf}\ref{stf2}, we get trees 
$S',T'\in\sct \jpi m$ satisfying  
$S'\nq m S$, $T'\nq m T$, $C=\raw{S'}s$, $D=\raw{T'}t$,  
so that $[\raw{T'}t]\cap(\ima f{[\raw{S'}s]})=\pu$.
Iterate this construction by exhaustion of 
all pairs $s,t\in2^m$. 
We get trees $S^\ast,T^\ast\in\sct \jpi m$ 
such that $S^\ast\nq m S$, $T^\ast\nq m T$, and 
$[\raw{T^\ast}t]\cap(\ima f{[\raw{S^\ast}s]})=\pu$ 
for all $s,t\in2^m$, that is  
$[{T^\ast}]\cap(\ima f{[S^\ast]})=\pu$.
Nodify $\vpi$ by  
$\xT\vpi Km=S^\ast$ instead of $S$ and 
$\xT\vpi Lm=T^\ast$ instead of $T$. 
Now \ref{ww43} holds for the modified $\vpi$.  

{\ubf  Case 2:} $K=L$, and then 
$S=T=\xT\vpi Km=\xT\vpi Lm$. 
It suffices to show that if  $s,t\in2^n$ 
(strings of length $n$) and $\raw Ss\sq T_0$, 
then there is a tree $S'\in \sct Pn$ such that  
$S'\nq m S$ and 
$[\raw{S'}t]\cap(\ima f[\raw{S'}s])=\pu$. 
Then the iteration by exhaustion of all 
those pairs of strings $s,t$ yields a tree   
$S^\ast\in\sct \jpi m$ such that  
$S^\ast\nq m S$ and  
$[\raw{T^\ast}t]\cap(\ima f{[\raw{S^\ast}s]})=\pu$ 
for all $s,t\in2^m$ with $\raw Ss\sq T_0$, that is, 
$[{S^\ast}]\cap(\ima f{[\raw{S^\ast}s]})=\pu$ 
whenevew $s\in2^m$ satisfies  
$\raw Ss\sq T_0$.
To achieve \ref{ww43}, it remains  
to modify the system $\vpi$ by  
$\xT\vpi Km=S^\ast$ instead of $S$. 

Thus let us carry out the construction of $S'$ 
for a pair of strings  
$s,t\in2^n$ with $\raw Ss\sq T_0$. 
It follows that $f$ 
is regular on $\raw Ss$ inside $\jpi$.
By Lemma~\ref{sadd}\ref{sadd0}, we have  
$\raw Ss=\req Su$ and $\raw St=\req Sv$, 
where $u,v\in T$ are strings of length 
$h=\lh u=\lh v$, and   
$\raw Ss=\tau\app(\raw St)$, 
where $\tau=u\app v$, by  
Lemma~\ref{trans}\ref{trans1}.
Lemma~\ref{rc}\ref{rc2} yields a pair of  
trees $U,V\in\jpi$ satisfying  
$U\sq \raw Ss$, $V\sq \raw St$, $V=\tau\app U$, and   
$[V]\cap(\ima f[U])=\pu$. 
Now, twice reducing the tree $S$ by means of  
Lemma~\ref{stf}\ref{stf2}, we get a tree  
$S'\in \sct Pn$ such that $S'\nq m S$ and  
$\raw{S'}s=U$, $\raw{S'}t=V$, so that  
$[\raw{S'}t]\cap(\ima f[\raw{S'}s])=\pu$, 
as required.  
\epF{Lemma}

Now return to the verification of \ref{dfex}\ref{dfex6}.
By the lemma, at least one  
\mus\ $\vpi(j)$  belongs to $\Da$, that is, conditions   
\ref{ww42}, \ref{ww41}, \ref{ww43} 
are satisfied for $\vpi=\vpi(j)$. 
The tree $V=\qfd{}L$ satisfies  
$V\sq T=\xT{\vpi(j)}Lm$. 
Moreover, the tree 
$U=\qfs{}K{s_0}$ satisfies 
$U\sq S=\bigcup_{s\in\Sg}\xT{\vpi(j)}Ks$, 
where $\Sg=\ens{s\in2^m}{\raw Ss\sq T_0}$, 
by \ref{ww42}, \ref{ww41}. 
However $[T]\cap(\ima f[S])=\pu$ 
by \ref{ww43}, 
thus $[V]\cap(\ima f[U])=\pu$, 
as required.
\vom
 
\qeD{Theorem~\ref{tex}}

\parf{The blocking sequence of \sfo s} 
\las{vop}

{\ubf We argue in the constructible universe $\rL$} 
in this section.

The forcing to prove~Theorem~\ref{Tun} 
will be defined as the union of a \dd\omi sequence  
of countable \sfo s, increasing in the sense of 
a relation $\prok$  
(Definition~\ref{dfex}). 
We here introduce the notational system to be used 
in this construction.

Let $\stf$ be the set of all \rit{countable} \sfo s. 
\index{forcing!$\stf$} 
\index{zzSTF@$\stf$} 
If $\vjpi=\sis{\jpi_\xi}{\xi<\omi}$ 
is a transfinite sequence of countable \sfo s, 
of length $\dom\vjpi=\la<\omi$, 
then let $\bigcup\vjpi=\bigcup_{\xi<\la}\jpi_\xi$, 
and let $\cM(\vjpi)$ be the least transitive model  
\index{zzMP@$\cM(\vjpi)$}%
\index{model!zzMP@$\cM(\vjpi)$}%
of $\zfcm$ of the form $\rL_\vt$, containing $\vjpi$, 
in which $\la$ and $\bigcup\vjpi$ are countable.

\bdf
\lam{vstf}
If $\al\in\Ord$ then let $\vstf_\al$ 
be the set of all \dd\la sequences    
$\vjpi=\sis{\jpi_\xi}{\xi<\al}$ 
of forcings $\jpi_\xi\in\stf$, 
satisfying the following:
\ben
\fenu
\itlb{vstf*} 
if\, $\ga<\dom\vjpi$\, then\, 
${\bigcup{(\vjpi\res\ga)}}\prol{\cM(\vjpi\res\ga)} 
\jpi_\ga$. 
\een
Let $\vstf=\bigcup_{\al<\omi}\vstf_\al$.
\edf

The set $\vstf\cup\vsto$ is ordered by the 
extension relations $\su$ and $\sq$. 

\ble
\lam{202}
Assume that\/ $\ka<\la<\omi$, 
and\/ $\vjpi=\sis{\jpi_\xi}{\xi<\ka}\in\vstf$. 
Then$:$
\ben
\renu
\itlb{2021}%
the union\/ $\jpi=\bigcup\vjpi$ is a countable\/ 
\sfo$;$

\itlb{2022}%
there is a sequence\/ $\vjqo\in\vstf$ 
such that\/ $\dom(\vjqo)=\la$ and\/   
$\vjpi\su\vjqo\;.$
\vyk{
\itlb{2023}%
если\/ $\vjpi,\vjqo,\vjR\in\vstf$ 
и\/ $\vjpi\su_\cM\vjqo\su\vjR$, то\/ 
$\vjpi\su_\cM\vjR\;.$
}
\een
\ele
\bpf
To prove \ref{2022} apply Theorem~\ref{tex} by induction 
on $\la$.
\epf

\bdf
[key definition]
\lam{dzap}
A sequence\/ $\vjpi\in\vstf$
\rit{blocks} a set $W\sq\vstf$ if either  
\index{blocks}%
\index{blocks!positively}%
\index{blocks!negatively}%
$\vjpi\in W$ (the positive block case)
or there is no sequence 
$\vjqo\in W$ satisfying $\vjpi\sq\vjqo$ 
(the negative block case).
\edf



Approaching the next 
\rit{blocking sequence theorem}, 
we recall that $\hc$ is the set of all 
\index{hereditarily countable!$\hc$}%
\index{zzhc@$\hc$}%
\rit{hereditarily countable\/} sets, so that 
$\hc=\rL_{\omi}$ in $\rL$.
See \cite[Part B, Chap.\,5,\;Section\,4]{skml} 
on definability classes $\is X n,\,\ip X n,\,\id X n$ 
\index{definability classes!$\is\hc n,\,\ip\hc n,\,\id\hc n$}%
\index{zzSHC@$\is\hc n,\,\ip\hc n,\,\id\hc n$}%
for any set $X$,
and especially on $\is\hc n,\,\ip\hc n,\,\id\hc n$ for 
$X=\hc$ in \cite[Sections 8, 9]{skmlD} or elsewhere.

\bte
[blocking sequence theorem, in $\rL$]
\lam{ep}
If\/ $\nn\ge3$ then there is a sequence\/ 
$\vdp=\sis{\dP_\xi}{\xi<\omi}\in\vsto$  
satisfying the following two conditions$:$
\ben
\renu
\itlb{ep1}%
$\vdp$, 
as the set of pairs\/ $\ang{\xi,\dP_\xi}$, 
belongs to the definability class\/ $\id\hc{\nn-1}\;;$

\itlb{ep2}%
{\rm(genericity of $\vjpi$ \poo\ $\is\hc{\nn-2}(\hc)$ sets)} 
\ if\/ $W\sq\vstf$ is a\/ $\is\hc{\nn-2}(\hc)$ set\/ 
{\rm(that is parameters from $\hc$ are admitted)}, 
\index{definability classes!$\is\hc n(\hc),\,\ip\hc n(\hc)$}%
\index{zzSHChc@$\is\hc n(\hc),\,\ip\hc n(\hc)$}%
then there is an ordinal\/ $\ga<\omi$ such that 
the restricted sequence\/ 
$\vdp\res\ga=\sis{\dP_\xi}{\xi<\ga}\in\vstf$
blocks\/ $W$.
\een
\ete
\bpf
Let $\lel$ be the canonical $\id{}1$ wellordering of $\rL$; 
thus its restriction to $\hc=\rL_{\omi}$ is $\id\hc1$.
As $\nn\ge3$, there exists a universal $\is\hc{\nn-2}$ set 
$\gU\sq\omi\ti\hc$. 
That is, $\gU$ is $\is\hc{\nn-2}$ 
(parameter-free $\is{}{\nn-2}$ definable in $\hc$), 
and for every set $X\sq\hc$ of class $\is\hc{\nn-2}(\hc)$ 
(definable in $\hc$ by a $\is{}{\nn-2}$ formula with arbitrary 
parameters in $\hc$) 
there is an ordinal $\al<\omi$ such that $X=\gU_\al$, 
where $\gU_\al=\ens{x}{\ang{\al,x}\in\gU}$). 
The choice of $\omi$ as the domain of parameters is 
validated by the assumption $\rV=\rL$, which implies the 
existence of a $\id\hc1$ surjection $\omi\onto\hc$.

Coming back to Definition~\ref{dzap}, note that for any 
sequence $\vjpi\in\vstf$ and any set $W\sq\vstf$ 
there is a sequence $\vjqo\in\vstf$ which satisfies 
$\vjpi\su\vjqo$ and blocks $W$.
This allows us to define $\vjqo_\al\in\vstf$ by induction 
on $\al<\omi$ so that $\vjqo_0=\pu$, 
$\vjqo_\la=\bigcup_{\al<\la}\vjqo_\al$, and each 
$\vjqo_{\al+1}$ is equal to the \dd\lel least 
sequence $\vjqo\in\vstf$ which satisfies 
$\vjqo_\al\su\vjqo$ and blocks $\gU_\al$. 
Then $\vdp=\bigcup_{\al<\omi}\vjqo_\al\in\vsto$. 

Condition \ref{ep2} holds by construction, while  
\ref{ep1} follows by a routine verification, 
based on the fact that $\vstf\in\id\hc1$.
\epf

\bdf
[in $\rL$]
\lam{vdp}
We fix a natural number $\nn\ge3$, for which  
\index{zzn@$\nn\ge 3$}%
\index{n@number $\nn\ge 3$}%
Theorem~\ref{Tun} is to be established.
We also fix a sequence  
$\vdp=\sis{\dP_\xi}{\xi<\omi}\in\vsto$, 
given by Theorem~\ref{ep} for this $\nn$.

If $\al<\omi$ then let $\cM_\al=\cM(\vdp\res\al)$ 
\imar{vdp, mdp al}
\index{zzMalpha@$\cM_\al$}%
\index{model!zzMalpha@$\cM_\al$}%
and $\mdp{\al}=\bigcup_{\xi<\al}\dP_\xi$. 

Let  
$\dP=\bigcup_{\xi<\omi}\dP_\xi$. 
\index{zzPd@$\dP$}%
\index{forcing!zzPd@$\dP$}%
\edf

\parf{CCC and some other forcing properties} 
\las{pres+}

The \sfo\ $\dP$ defined by \ref{vdp} 
will be the forcing notion for the proof 
of Theorem~\ref{Tun}. 
Here we establish some forcing properties of  
$\dP$, including CCC.

{\ubf We continue to argue in the conditions and 
notation of Definition~\ref{vdp}}.

\ble
\lam{jden}
$\dP$ and all sets\/ $\dP_\xi$, 
$\mdp\al$ are\/ \sfo s. 
In addition$:$
\ben
\renu
\itlb{prop1}%
if\/ $\al<\omi$ then\/ ${\mdp\ga}\prol{\cM_\ga}\dP_\ga\,;$

\itlb{jden1}%
if\/ $\al<\omi$ and the set\/ 
\imar{jden}
$D\in\cM_\al\yt D\sq\mdp\al$ 
is pre-dense in\/ $\mdp\al$ then it is 
pre-dense in\/ $\dP$, too$;$ 

\itlb{jden2}%
every set\/ $\dP_\al$ is pre-dense in\/ $\dP\,;$


\itlb{prop3}%
if\/ $Q\sq\pes$ belongs to\/ 
$\is\hc{\nn-2}(\hc)$ and\/ 
$Q^-=\ens{T\in\pes}{\neg\:\sus S\in Q\,(S\sq T)}$  
then\/ $\dP\cap{(Q\cup Q^-)}$ is dense in\/ $\dP$$;$

\itlb{prop4}%
if\/ $\rc\in \ccf$ is a code of continuous function 
and\/ 
$$
\cb\rc=\ens{T\in\pes}
{f_\rc\res[T]\text{\rm\ is a constant or a bijection}}
$$
then the set\/ $\dP\cap{\cb\rc}$ is dense in\/ $\dP$$;$

\itlb{prop5}%
let\/ $D=\ens{T\in\pes}{\oi T\text{ is co-infinite}}$  
{\rm(see Definition~\ref{std} on $\oi T$),}  
then 
the set\/ $\dP\cap D$ is dense in\/ $\dP$.
\een
\ele
\bpf
\ref{prop1}
holds by \ref{vstf*} of Definition~\ref{vstf}.

\ref{jden1}  
We use induction on $\ga\yd \al\le\ga<\omi$, 
to check that if $D$ is pre-dense in $\mdp\ga$  
then it remains pre-dense in   
$\mdp\ga\cup\dP_\ga=\mdp{\ga+1}$ by \ref{prop1}
and \ref{dfex}\ref{dfex3}. 
Limit steps, including the final step to $\dP$ 
($\ga=\omi$) are routine. 

\ref{jden2}
$\dP_\al$ is dense in  
$\mdp{\al+1}=\mdp\al\cup\dP_\al$ 
by \ref{dfex}\ref{dfex1}.
It remains to refer to \ref{jden1}.


\ref{prop3}
Let $T_0\in\dP$, that is, 
$T_0\in\mdp{\al_0}$, $\al_0<\omi$. 
The set $W$ of all sequences $\vjpi\in\vstf$, 
such that $\vdp\res\al_0\sq\vjpi$ and 
$\sus T\in Q\cap(\bigcup\vjpi)\,(T\sq T_0)$, 
belongs to $\is\hc{\nn-2}(\hc)$ 
along with $Q$. 
Therefore there is an ordinal $\al<\omi$ such that 
$\vdp\res\al$ blocks $W$. 
We have two cases.

{\ubf Case 1:\/} 
$\vdp\res\al\in W$.
Then the related tree $T\sq T_0$ belongs to 
$Q\cap\dP$.

{\ubf Case 2:\/} 
there is no sequence in $W$ which extends  
$\vdp\res\al$.
Let $\ga=\tmax\ans{\al,\al_0}$. 
Then $\mdp\ga\prol{\cM_\ga}\dP_\ga$ by \ref{prop1}. 
As $\al_0<\ga$, there is a tree $T\in\dP_\ga$, 
$T\sq T_0$. 
We claim that $T\in Q^-$, which completes the 
proof in Case 2. 

Suppose to the contrary that $T\nin Q^-$, 
thus there is a tree $S\in Q$, $S\sq T$.
The set 
$\bR=\dP_\ga\cup\ens{\sg\app (\raw St)}{t,\sg\in\bse}$ 
is a countable \sfo\ and   
$\mdp\ga\prol{\cM_\ga}\bR$ by Lemma~\ref{tec}\ref{tec3}. 
It follows that the sequence $\vjR$ defined by  
$\dom\vjR=\ga+1$, 
$\vjR\res\ga=\vdp\res\ga$, and $\vjR(\ga)=\bR$, 
belongs to $\vstf$, and even $\vjR\in W$ since 
$S\in Q\cap\bR$.  
Yet $\vdp\res\al\prol{\cM_\ga}\vjR$ by construction, 
which contradicts to the Case 2 hypothesis. 

\ref{prop4}
A routine estimation gives $\cb\rc\in\is\hc1(\ans{\rc})$. 
The set $\cb\rc$ is dense in $\pes$
by Corollary~\ref{ssq}, 
thus $(\cb\rc)^-$ is empty.
Now the result follows from \ref{prop3}.

\ref{prop5}
Similarly to \ref{prop4}, $D\in\is\hc1$   
and $D$ is dense in $\pes$.
\epf

\bcor
\lam{inc}
If\/ $\al<\omi$ and trees\/ $T,T'\in\mdp{\al}$   
are incompatible in\/ $\mdp{\al}$ then\/ $T,T'$  
are incompatible in\/ $\dP$, too. 
\ecor
\bpf
Prove by induction on $\ga$ that if $\al<\ga\le\omi$ 
then $T,T'$ are incompatible in $\mdp{\ga}$, using 
Lemma~\ref{jden}\ref{prop1}, 
and Lemma~\ref{tec}\ref{tec4} on limit steps. 
\epf

To prove CCC we'll need the following lemma.

\ble
[in $\rL$]
\lam{club}
If\/ $X\sq\hc=\rL_{\omi}$ then the set\/ $\skri O_X$ 
of all ordinals\/ $\al<\omi$ such that the model\/ 
$\stk{\rL_\al}{X\cap\rL_\al}$  
is an elementary submodel of\/ $\stk{\rL_{\omi}}{X}$ 
and\/ $X\cap\rL_\al\in\cM_\al$, is unbounded in\/ $\omi$.
Generally if\/ $X_n\sq\hc$ for all\/ $n$  
then the set\/ $\skri O$ of all ordinals\/ 
$\al<\omi$ such that\/  
$\stk{\rL_\al}{\sis{X_n\cap\rL_\al}{n<\om}}$ is  
an elementary submodel of\/  
$\stk{\rL_{\omi}}{\sis{X_n}{n<\om}}$ 
and\/ $\sis{X_n\cap\rL_\al}{n<\om}\in\cM_\al$, 
is unbounded in\/ $\omi$.
\ele
\bpf
Let $\al_0<\omi$. 
There is a countable elementary submodel  $M$ 
of $\stk{\rL_{\om_2}}{\in}$  
which contains $\al_0\yi\omi\yi X$  
and is such that the set $M\cap\rL_{\omi}$ is transitive. 
Consider the Mostowski collapse $\phi:M\onto\rL_\la$.  
Let $\al=\phi(\omi)$. 
Then $\al_0<\al<\la<\omi$ and $\phi(X)=X\cap\rL_\al$ 
by the choice of $M$. 
We conclude that $\stk{\rL_\al}{X\cap\rL_\al}$ is  
an elementary submodel of $\stk{\rL_{\omi}}{X}$.
And $\al$ is uncountable in $\rL_\la$, hence  
$\rL_\la\sq\cM_\al$. 
It follows that $X\cap\rL_\al\in\cM_\al$, as   
$X\cap\rL_\al\in\rL_\la$ by construction.

The more general claim is proved similarly.
\epf

\bcor
[in $\rL$]
\lam{ccc}
$\dP$ is a CCC forcing, therefore\/ 
\dd\dP generic extensions preserve cardinals.
\ecor
\bpf
Consider any maximal antichain $A\sq\dP$. 
By Lemma~\ref{club} there is an orrdinal $\al$ 
such that $\stk{\rL_\al}{\dP',A'}$ is  
an elementary submodel of $\stk{\rL_{\omi}}{\dP,A}$, 
where $\dP'=\dP\cap\rL_\al$ and $A'=A\cap\mdp\al$, 
and in addition $\dP',A'\in\cM_\al$. 
By the elementarity, we have  
$\dP'=\mdp\al$ and $A'=A\cap\mdp\al\in\cM_\al$, 
and $A'$ is a maximal antichain, hence a pre-dense  
set, in $\mdp\al$. 
But then $A'$ is a pre-dense set, hence, 
a maximal antichain, in the whole set $\dP$ 
by Lemma~\ref{jden}\ref{jden1}. 
Thus $A=A'$ is countable.
\epf

\parf{Generic model} 
\las{bex1}

This section presents some properties of  
\dd\dP generic extensions $\rL[G]$ of $\rL$ 
obtained by adjoining a \dd\dP generic  
set $G\sq\dP$ to $\rL$. 
Recall that the forcing notion $\dP\in\rL$ 
was introduced by Definition~\ref{vdp}, along with 
some related notation. 

\vyk{
\bcor
\lam{gga}
If\/ $\al<\omil$ and a set\/  
$G\sq\dP$ is\/ \dd\dP generic over\/ $\rL$  
then the set\/ $G'=G\cap\mdp\al$ is\/ 
\dd{\mdp\al}generic over\/ $\cM_\al$.
\ecor
\bpf
The trees in $G'$ are pairwise compatible inside 
$\mdp\al$ by Corollary~\ref{inc}. 
Moreover if a set $D\in\cM_\al$, $D\sq\mdp\al$ 
is dense in $\mdp\al$ then it is pre-dense in $\dP$
by Lemma~\ref{jden}. 
It follows that $G\cap D\ne\pu$, 
therefore $G'\cap D\ne\pu$.
\epf
}

The next lemma involves the coding system for 
continuous maps introduced in Section~\ref{nef}. 
If $G\sq\dP$ is generic over $\rL$ 
and $\rc\in\ccf\cap\rL$, 
then define $\rc[G]=f_\rc(\xx G)\in\dn\cap\rL[G]$; 
the definition of $\xx G$ see \ref{xg}.

\ble
[continuous reading of names]
\lam{repd}
If a set $G\sq\dP$ is generic over\/ $\rL$ 
and\/ $x\in\dn\cap\rL[G]$ then 
there exists a code\/ $\rc\in\ccf\cap\rL$  
such that\/ $x=\rc[G]$.
\ele
\bpf
One of basic forcing lemmas 
(Lemma 2.5 in \cite[Chap.\,4]{skml})  
claims that there is a \dd\dP name $t\in\rL$ for $x$, 
satisfying $x=t[G]$ 
(the \dd Gvaluation of $t$), 
and it can be assumed that $\dP$ forces that $t$   
is valuated as a real in $\dn$. 
Then the sets 
$F_{ni}=\ens{T\in\dP}{T\text{ forces }t(n)=i}$ 
($n<\om$ and $i=0,1$) 
satisfy the following:
\ben
\nenu
\itla{FF1}
the indexed set $\sis{F_{ni}}{n<\om\land i=0,1}$ 
belongs to $\rL$;

\itla{FF2}
if $n<\om$, $S\in F_{n0}$, $T\in F_{n1}$, then $S,T$ 
are incompatible in $\dP$;

\itla{FF3}
if $n<\om$ then the set $F_n= F_{n0}\cup F_{n1}$ 
is open dense in $\dP$.
\een

{\ubf We argue in\/ $\rL$}.
Pick a maximal antichain $A_n\sq F_n$ in each $F_n$. 
Then all sets $A_n$ are maximal antichains in $\dP$ 
by \ref{FF3}, and all $A_n$ are countable by 
Corollary~\ref{ccc}. 
Therefore there is an ordinal $\al<\omil$ such that 
the set $\bigcup_nA_n\sq\mdp\al$
and the sequence $\sis{A_{n}}{n<\om}$ belong to $\cM_\al$. 
Note that $G\cap\dP_\al\ne\pu$ by 
Lemma~\ref{jden}\ref{jden2}; let $U\in {G\cap\dP_\al}$. 
As ${\mdp\al}\prol{\cM_\al}\dP_\al$ 
by Lemma~\ref{jden}\ref{prop1}, 
we have $U\sqf \bigcup A_n$ for all $n$, there exists 
a finite set  
$A'_n\sq A_n$ such that $U\sq \bigcup A'_n$.

Let $A'_{ni}=A'_{n}\cap F_{ni}$ and 
$X_{ni}=[U]\cap \bigcup_{T\in A'_{ni}}[T]$, $i=0,1$. 
We claim that $X_{n0}\cap X_{n1}=\pu$. 
Indeed otherwise there exist trees $T_i\in F_{ni}$ 
such that $Z=[U]\cap[T_0]\cap [T_1]$ is non-empty.
By Lemma~\ref{tec}\ref{tec1}, $Z$ is clopen in $[U]$. 
Therefore there is a tree $U'\in\dP_\al$ such that 
$[U']\sq Z$, hence, $T_0$ and $T_1$ 
are compatible in $\dP$, 
which contradicts \ref{FF2} by construction. 
Thus indeed $X_{n0}\cap X_{n1}=\pu$. 

As clearly $X_{n0}\cup X_{n1}=[U]$, we can define 
a continuous $g:[U]\to\dn$ such that  
$g(x)(n)=i$ iff $x\in X_{ni}$.
The map $g$ can be extended to a continuous  
$f:\dn\to\dn$, that is $f(x)(n)=i$ whenever 
$x\in X_{ni}$ --- for all $n<\om$ and $i=0,1$.
We have $f=f_\rc$, where $\rc\in\ccf\cap\rL$.

{\ubf We argue in\/ $\rL[G]$}.
We skip a routine verification of $x=\rc[G]$.
\epf

\ble
\lam{sym}
If\/ $G\sq\dP$ is generic over\/ $\rL$ then\/ $\xx G$ 
is not\/ 
\od\ in\/ $\rL[G]$.
\ele
\bpf
Assume to the contrary that 
$\vt(x)$ is a formula with ordinal parameters,  
and a tree $T\in G$ \dd\dP forces that  
$\xx G$ is the only real $x \in \dn$ 
satisfying $\vt(x)$. 
Let $s =\roo T$ and $n = \lh s$. 
Then $T$ contains both $s\we 0$ and $s\we1$. 
Then either $s\we0\su \xx G$ or $s\we1\su \xx G$. 
Let say $s\we0\su \xx G$. 

Let $\sg = 0^n\we1$, so that the strings 
$s\we0$, $s\we1$, $\sg$ belong to $2^{n+1},$   
$s\we 1 = \sg\app s\we 0$, and $\sg\app T=T$ 
by Lemma~\ref{trans}\ref{trans1}. 
As $\dP$ is invariant under the action of $\sg$,
the set $G' = \sg\app G=\ens{\sg\app S}{S\in G}$ 
is \dd\dP generic over $\rL$, and $T = \sg\app T\in G'$. 
We conctude that it is true in $\rL[G'] = \rL[G]$ that   
$\xx{G'}= \sg\app \xx G$ is still the unique real in 
$\dn$ satisfying $\vt(\xx{G'})$. 
But $\xx{G'}\ne \xx{G}$!
\epf

\parf{Definability of the set of generic reals} 
\las{bex2}

We continue to argue in the context of 
Definition~\ref{vdp}. 
The goal of this section is to study the definability 
of the set of all \dd{\dP}generic reals $x\in\dn$ 
in \dd\dP generic extensions of $\rL$.
 
\ble
\lam{mod1}
In an transitive model of\/ $\zf$ extending\/ $\rL$, 
it is true that a real\/ $x\in\dn$ is\/ 
\dd{\dP}generic over\/ $\rL$ iff\/ 
$x$ belongs to the set\/ $\gen_\dP=
\bigcap_{\al<\omil}\bigcup_{T\in\dP_{\al}}[T]$.%
\ele
\bpf
All sets $\dP_{\al}$ are pre-dense in  
$\dP$ by Lemma~\ref{jden}\ref{jden2}.
Therefore all \dd{\dP}generic reals belong   
to $\gen_\dP$. 
On the other hand, any maximal antichain  
$A\in\rL\yd A\sq\dP$ is countable in $\rL$ 
by Corollary~\ref{ccc}, 
and hence $A\sq\mdp\al$ and $A\in\cM_\al$ 
for some index $\al<\omil$.
But then every tree $T\in\dP_{\al}$ satisfies  
$T\sqf\bigcup A$ by Lemma~\ref{tec}\ref{tec2}. 
We conclude that  
$\bigcup_{T\in\dP_{\al}}[T]\sq \bigcup_{S\in A}[S]$.
\epf       

According to the next lemma, \dd\dP generic 
extensions do not contain
\dd{\dP}generic reals, 
except the real $\xx G$ itself and reals connected 
to $\xx G$ in terms of the equivalence relation $\Eo$ 
(see Footnote~\ref{deo}).
We observe that the \dd\Eo class 
$$
[x]_{\Eo}=\ens{y\in\dn}{x\Eo y}=
\ens{y\in\dn}{\sus s\in\bse(y=s\app x)}
$$ 
\index{equivalence relation!$[x]_{\Eo}$, \dd\Eo class}%
\index{zzxE0@$[x]_{\Eo}$}%
of any real $x\in\dn$ is a countable set.

\ble
\lam{only}
Let\/ $G\sq\dP$ be a\/ \dd\dP generic set  
over\/ $\rL$. 
Then it is true in\/ $\rL[G]$ that\/ $\gen_\dP=\eko{\xx G }$.
\ele
\bpf
The real $\xx G $ is \dd{\dP}generic, 
hence $\xx G\in\gen_\dP$ by Lemma~\ref{mod1}. 
Yet every set $\dP_\al$ is a \sfo, 
that is by definition it is closed under the action 
$s\app T$ of any string $s\in\bse.$ 
This implies $\eko{\xx G}\sq\gen_\dP$.

To prove in the other direction, assume to the  
contrary that  $x\in\rL[G]\cap\dn$, 
$x\in \gen_\dP\bez\eko{\xx G}$.
By Lemma~\ref{repd}, 
there is a code $\rc\in\rL\cap\ccf$ such that 
$x=\rc[G]$. 
By the contrary assumption there is a tree $T_0\in G$ 
which forces
\bce 
$\rc[\uG]\in\gen_\dP\bez\eko{\xx{\uG}}$.
\ece  
Fix an ordinal $\al<\omil$ such that $\rc\in\cM_\al$. 
We claim that (in $\rL$) 
\rit{the function\/ $f=f_\rc$ 
is regular on\/ $T_0$ inside\/ $\mdp\al$}.

Indeed otherwise there exist 
$\sg\in\bse$ and $T\in\mdp\al$ such that $T\sq T_0$  
and $\sg\app f(x)=x$ for all $x\in T$. 
Then $T$ forces $\sg\app \rc[\uG]=\xx{\uG}$, 
that is, forces $\rc[\uG]\in\eko{\xx{\uG}}$, 
contrary to the choice of $T_0$. 
The regularity is established. 

Recall that $\mdp\al\prol{\cM_\al}\dP_\al$ by 
Lemma~\ref{jden}\ref{prop1}. 
Therefore by \ref{dfex}\ref{dfex1},  
there is a tree $U\in\dP_\al$, $U\sq T_0$.
And by \ref{dfex}\ref{dfex6} if $V\in\dP_\al$  
then $[V]\cap (\ima f{[U]})=\pu$. 
Thus if $V\in\dP_\al$ then $U$ forces 
$\rc[\uG]\nin [V]$, hence forces  
$\rc[\uG]\nin \gen_\dP$,  
which contradicts to the choice of $T_0$. 
\epf

\bcor
\lam{gdef}
Let\/ $G\sq\dP$ be a \dd\dP generic set  
over\/ $\rL$. 
Then the\/ \dd{\Eo}class\/ $\eko{\xx G}$ is a\/ 
$\ip1\nn$ set in\/ $\rL[G]$.
\ecor
\bpf
A routine verification of $\gen_\dP\in\ip1\nn$ in $\rL[G]$ 
using the property \ref{ep}\ref{ep1} of the sequence 
$\vjpi$ is left to the reader.
\epf

\parf{Auxiliary forcing relation} 
\las{fs}


Here we introduce a key tool for the proof of 
claim \ref{Tun1} of Theorem~\ref{Tun}. 
This is a forcing-like relation $\fo$. 
It is not explicitly connected with the forcing notion 
$\dP$ 
(but rather connected with the full 
Silver forcing $\pes$), 
however it will be compatible with $\dP$ for  
formulas of certain quantifier complexity 
(Lemma~\ref{ft.}). 
The crucial advantage of $\fo$ will be its invariance  
under two certain groups of transformations  
(Lemma~\ref{inv}), 
a property that cannot be expected for $\dP$.
This will be the key argument in the proof 
of Theorem~\ref{lun} below.

{\ubf We argue in $\rL$.}

We consider a \rit{language\/} $\cL$ containing  
\index{formula!language $\cL$}%
variables $i,j,k,\dots$ of type 0 with the domain $\om$ 
and variables $x,y,z,\dots$ of type 1 with the domain $\dn$. 
\rit{Terms} are variables of type 0 and    
expressions like $x(k)$. 
Atomic formulas are those of the form $R(t_1,\dots,t_n)$, 
where $R\sq\om^n$ is any \dd nary relation  
on $\om$ in $\rL$. 
A formula is \rit{arithmetic} if it does not contain 
\index{formula!arithmetic}%
variables of type 1. 
Formulas of the form
$$
\sus x_1\:\kaz x_2\:\sus x_3\:\dots\:\sus(\kaz)\,x_n\,\Psi
\quad\text{and}
\quad
\kaz x_1\:\sus x_2\:\kaz x_3\:\dots\:\kaz(\sus)\,x_n\,\Psi\,,
$$
where $\Psi$ is arithmetic, are of types  
\index{formula!$\lis1n$, $\lip1n$}%
\index{zzLS1n@$\lis1n$}%
\index{zzLP1n@$\lip1n$}%
$\lis1n$, resp., $\lip1n$. 

In addition we allow codes  
$\rc\in\ccf$ to substitute free variables of type 1. 
The semantics is as follows. 
Let $\vpi:=\vpi(\rc_1,\dots,\rc_k)$ be an \dd\cL formula, 
with all codes in $\ccf$ explicitly indicated, and   
let $x\in\dn$. 
\index{formula!$\vpi[x]$}%
\index{zzfix@$\vpi[x]$}%
By $\vpi[x]$ we denote the formula 
$\vpi(f_{\rc_1}(x),\dots,f_{\rc_k}(x))$, 
where all $f_{\rc_i}(x)$ are reals in $\dn$, of course. 
 
Arithmetic formulas and those in $\lis1n\cup\lip1n$, 
$n\ge1$, will be called \rit{normal}. 
\index{formula!normal}%
If $\vpi$ is a formula in $\lis1n$ or $\lip1n$ then     
$\vpi^-$ is the result of canonical transformation of 
\index{formula!$\vpi^-$}%
\index{zzfi@$\vpi^-$}%
$\neg\:\vpi$ to $\lip1n$, resp., $\lis1n$ form. 
For arithmetic formulas, let $\vpi^-:=\neg\:\vpi$.

\bdf
[in $\rL$]
\lam{d:fo}
We define a relation $T\fo\vpi$ between trees $T\in\pes$ 
and closed normal \dd\cL formulas:
\ben
\Renu
\itlb{fo1}%
if $\vpi$ is a closed \dd\cL formula, arithmetic or 
in $\lis11\cup\lip11$, 
then $T\fo\vpi$ iff $\vpi[x]$ holds for all  
$x\in[T]$;

\itlb{fo2}%
if $\vpi:=\sus x\,\psi(x)$ is a closed  
$\lis1{n+1}$ formula, $n\ge 1$ 
($\psi$ being of type $\lip1n$),
then $T\fo\vpi$ iff there is a code $\rc\in\ccf$   
such that $T\fo\psi(\rc)$;

\itlb{fo3}%
if $\vpi$ is a closed $\lip1{n}$ formula, $n\ge 2$,
then $T\fo\vpi$ iff there is no tree $S\in\pes$   
such that $S\sq T$ and $S\fo\psi^-$.
\een
Let $\Fo\vpi=\ens{T\in\pes}{T\fo\vpi}$ and 
$\des\vpi=\Fo\vpi\cup\Fo{\vpi^-}$.
\edf

\ble
[in $\rL$]
\lam{deff}
If\/ $m\ge2$ and\/ $\vpi$ is a closed formula in\/ 
$\lis1m$, resp., $\lip1m$, then the set\/ 
$\Fo\vpi$ belongs to\/ $\is\hc{m-1}(\hc)$, 
resp., $\ip\hc{m-1}(\hc)$.
\ele
\bpf
For $\lip11$ formulas, Definition~\ref{d:fo}\ref{fo1}
implies $\Fo\vpi\in\fp11$, 
thus $\Fo\vpi$ belongs to $\id\hc{1}(\hc)$. 
Then argue by 
induction on Definition~\ref{d:fo}\ref{fo2},\ref{fo3}.
\epf

Recall that a number $\nn\ge3$ is fixed by  
Definition~\ref{vdp}.

\ble
[in $\rL$]
\lam{dens}
Let\/ $\vpi$ be a closed normal\/ \dd\cL formula. 
Then the set\/ $\des\vpi$ is dense in\/ $\pes$. 
If\/ $\vpi$ is of type\/ $\lis1m$, 
$m<\nn$, then\/ 
$\des\vpi\cap\dP$ is dense in\/ $\dP$.
\ele
\bpf
It suffices to establish the density of $\des\vpi$ 
for formulas $\vpi$ as in \ref{fo1}. 
If $\vpi$ is such then the set   
$X(\vpi)=\ens{x\in T}{\vpi[x]}$ 
belongs to $\fs11\cup\fp11$, that is, it has the Baire  
property inside $[T]$. 
Therefore at least one of the two complimentary sets  
$X(\vpi)$, $X(\vpi^-)$ is not meager in $[T]$. 
It remains to apply Corollary~\ref{ssq}. 

The second claim follows from the first one by lemmas  
\ref{deff} and \ref{jden}\ref{prop3}.
\epf

\vyk{
\ble
[in $\rL$]
\lam{lres}
Let\/ $\vpi:=\vpi(\rc_1,\dots,\rc_m)$ be 
a closed normal\/ \dd\cL formula, 
with all codes in\/ $\ccf$ explicitly indicated,
and\/ $\vpi':=\vpi(\rc'_1,\dots,\rc'_m)$, where\/ 
$\rc'_1,\dots,\rc'_m$ is another string of codes  
in\/ $\ccf$.
Assume that\/ $T\in\pes$ and\/ $\vpi[x]$ 
coincides with\/ $\vpi'[x]$ for all\/ $x\in[T]$.
Then\/ $T\fo\vpi$ iff\/ $T\fo\vpi'$.
\ele


\bpf
An elementary proof by induction on the 
complexity of $\vpi$ is left to the reader.
\epf
}

\parf{Invariance} 
\las{inva}

It happens that the relation $\fo$ is invariant under 
two rather natural groups of transformations of $\pes$. 
Here we prove the invariance.
{\ubf We still argue in $\rL$.}

{\ubf First group.} 
Let $h\sq\om$. 
If $x\in\dn$ then a real $h\app x\in\dn$ is defined 
\index{action!$h\app x$}%
by $(h\app x)(j)=1-x(j)$ for $j\in h$, but 
$(h\app x)(j)=x(j)$ for $j\nin h$. 
If $X\sq\dn$ then let $h\app X=\ens{h\app x}{x\in X}$. 
\index{action!$h\app X$}%
Accordingly, if $s\in\bse$ and $n=\lh s$ then a string 
$h\app s\in\bse$ is defined by  
\index{action!$h\app s$}%
$\dom(h\app s)=n=\dom s$ and if $j<n$ then 
$(h\app x)(j)=1-x(j)$ for $j\in h$, but  
$(h\app x)(j)=x(j)$ for $j\nin h$.
If $T\sq\bse$ then let $h\app T=\ens{h\app s}{s\in T}$. 
\index{action!$h\app T$}%
Then obviously $T\in\pes$ iff  $h\app T\in\pes$. 

If $f:\dn\to\dn$ then a function $h\app f:\dn\to\dn$ 
is defined by $(h\app f)(x)=f(h\app x)$, 
equivalently, $(h\app f)(h\app x)=f(x)$.
If $f$ is continuous, then $f=f_\rc$, where  
$\rc\in\ccf$, and there is a canonical definition of  
a code $h\app\rc\in\ccf$ such that  
\index{action!$h\app \rc$}%
$h\app (f_\rc)=f_{h\app\rc}$.

Finally if $\vpi:=\vpi(\rc_1,\dots,\rc_k)$ is a  
\dd\cL formula then let $h\vpi$ be the formula 
\index{action!$h\vpi$}%
$\vpi(h\app\rc_1,\dots,h\app\rc_k)$. 
Then  
$(h\vpi)[h\app x]$ coincides with $\vpi[x]$. 

{\ubf Second group.} 
Let $\ib$ be the set of all idempotent 
bijections $b:\om\onto\om$, that is, we require 
that $b(j)=b\obr(j)$, $\kaz j$. 
If $x\in\dn$ then define $b\app x\in\dn$ by 
$(b\app x)(j)=x(b(j))$, $\kaz j$. 
Let $b\app X=\ens{b\app x}{x\in X}$, for $X\sq\dn$. 
If $T\in\pes$ then put 
$b\app T=\ens{x\res m}{x\in(b\app[T])\land m<\om}$.
Clearly $T\in\pes$ iff  $b\app T\in\pes$. 

If $f:\dn\to\dn$ then a function $b\app f:\dn\to\dn$ 
is defined similarly to the above 
by $(b\app f)(x)=f(b\app x)$, 
equivalently, $(b\app f)(b\app x)=f(x)$.
If $f$ is continuous, then $f=f_\rc$, where  
$\rc\in\ccf$, and still there is a canonical definition 
of a code $b\app\rc\in\ccf$ such that  
\index{action!$b\app \rc$}%
$b\app (f_\rc)=f_{b\app\rc}$.

And finally if $\vpi:=\vpi(\rc_1,\dots,\rc_k)$ is a  
\dd\cL formula then let $b\vpi$ be the formula 
\index{action!$b\vpi$}%
$\vpi(b\app\rc_1,\dots,b\app\rc_k)$. 
Then  
$(b\vpi)[b\app x]$ coincides with $\vpi[x]$.

\ble
[in $\rL$]
\lam{inv}
Let\/ $T\in\pes$ and\/ 
$\vpi$ be a closed normal\/ \dd\cL formula.
Then\/ 
\ben
\renu
\itlb{inv1}
if\/ $h\sq\om$ then\/  
$T\fo\vpi$ iff\/ ${(h\app T)}\fo h\vpi$$;$

\itlb{inv2}
if\/ $b\in\ib$ then\/  
$T\fo\vpi$ iff\/ ${(b\app T)}\fo b\vpi$.
\een
\ele
\bpf
\ref{inv1}
If $\vpi$ is of type \ref{d:fo}\ref{fo1} then it 
suffices to note that, first, 
$[h\app T]=\ens{h\app x}{x\in[T]}$, 
and second, if $x\in[T]$ then $\vpi[x]$ coincides with 
$(h\vpi)[h\app x]$. 
A routine induction based on   
Definition~\ref{d:fo}\ref{fo2},\ref{fo3} 
is left to the reader.
\epf

\bcor
\lam{tt'}
Let\/  
$\vpi$ be a closed normal\/ \dd\cL formula, such 
that if a code\/ $\rc$ occurs in\/ $\vpi$ then\/ 
$f_\rc$ is a constant.
Assume that\/ $S,T\in\pes$ and the splitting sets\/
$\oi S$, $\oi T$ {\rm(see Section~\ref{tre})} 
are both co-infinite.
Then\/ $S\fo\vpi$ iff\/ $T\fo\vpi$.
\ecor
\bpf
Assume that $S\fo\vpi$.
As both $\oi S$ and $\oi T$ are co-infinite 
(and they are infinite anyway), there is a bijection 
$b\in\ib$ such that $\ima b{(\oi S)}=\oi T$. 
Then the tree $S'=b\app S$ satisfies 
$\oi{S'}=\oi T$, 
and still $S'\fo\vpi$ by Lemma~\ref{inv}\ref{inv2}, 
as all 
occurring codes define constant functions, and hence 
$b\vpi$ and $\vpi$ essentially coincide. 
Now, as $\oi{S'}=\oi T$, 
there is a set $h\sq\om\bez\oi{S'}=\om\bez\oi T$
such that $T=h\app S'$. 
Then $T\fo\vpi$ by Lemma~\ref{inv}\ref{inv1} and the 
constantness of the codes involved, as above. 
\epf

\vyk{
\parf{Compatibility with the true forcing} 
\las{ctf}

We now show that the relation $\fo$ is compatible 
with the forcing notion $\dP$ with respect to 
formulas of type $\lis1\nn$ and lower.

\ble
\lam{ft.}
Assume that\/ $1\le n<\nn$, $\vpi\in\rL$ is a    
closed formula in\/ $\lip1{n}\cup\lis1{n+1}$, and a set\/ 
$G\sq\dP$ is generic over\/ $\rL$. 
Then the sentence\/ $\vpi[\xx G]$ is true in\/ $\rL[G]$  
if and only if\/ $\sus T\in G\,(T\fo\vpi)$.
\ele
\bpf
{\ubf Base of induction:} 
$\vpi$ is arithmetic or belongs to $\lis11\cup\lip11$,
as in \ref{d:fo}\ref{fo1}. 
If $T\in G$ and $T\fo\vpi$, then $\vpi[\xx G]$ 
holds by the Shoenfield absoluteness theorem, as 
$\xx G\in[T]$. 
The inverse holds by Lemma~\ref{dens}. 

{\ubf Step $\lip1n\imp\lis1{n+1}$.} 
Let $\vpi$ be $\sus x\,\psi(x)$ where $\psi$ 
is of type $\lip1n$. 
Assume that $T\in G$ and $T\fo\vpi$. 
Then by Definition~\ref{d:fo}\ref{fo2}  
there is a code $\rc\in\ccf\cap\rL$ such that 
$T\fo\psi(\rc)$.
By the inductive hypothesis, 
the formula $\psi(\rc)[\xx G]$, that is,
$\psi[G](f_\rc(\xx G))$, is true in $\rL[G]$. 
But then $\vpi[\xx G]$ is obviously true as well. 

Conversely assume that $\vpi[\xx G]$ is true. 
Then there is a real $y\in\rL[G]\cap\dn$ such that
$\psi[G](y)$ is true.
By Lemma~\ref{repd},  
$y=f_\rc(\xx G)$ for a code $\rc\in\ccf\cap\rL$. 
But then $\psi(\rc)[\xx G]$ is true in $\rL[G]$. 
By the inductive hypothesis, there is a tree $T\in G$  
satisfying $T\fo\psi(\rc)$. 
Then $T\fo\vpi$ as well. 

{\ubf Step $\lis1n\imp\lip1{n}$.} 
Let $\vpi$ be a $\lip1n$ formula, $n\ge2$. 
By Lemma~\ref{dens}, there is a tree $T\in G$  
such that either $T\fo\vpi$ or $T\fo\vpi^-$.
If $T\fo\vpi^-$ then $\vpi^-[\xx G]$ is true by  
the inductive hypothesis, hence $\vpi[\xx G]$ is false. 
Now assume that $T\fo\vpi$. 
We have to prove that $\vpi[\xx G]$ is true. 
Suppose otherwise. 
Then $\vpi^-[\xx G]$ is true. 
By the inductive hypothesis, there is a tree  
$S\in G$ such that $S\fo\vpi^-$.
But the trees $S,T$ belong to the same generic set 
$G$, hence they are compatible, which leads to a 
contradiction with the assumption $T\fo\vpi$, 
according to Definition~\ref{d:fo}\ref{fo3}.
\epf
}

\parf{The final argument} 
\las{zav}

Recall that $\nn\ge3$ is fixed by Definition~\ref{vdp}. 

The last part of the proof of Theorem~\ref{Tun} 
will be Lemma~\ref{lun}. 
Note the key ingredient of the proof: 
we surprisingly approximate the 
forcing $\dP$, definitely non-invariant under 
the transformations considered in Section~\ref{inva}, 
by the invariant relation $\fo$, using 
the next lemma.  

\ble
\lam{ft.}
Assume that\/ $1\le n<\nn$, $\vpi\in\rL$ is a    
closed formula in\/ $\lip1{n}\cup\lis1{n+1}$, and a set\/ 
$G\sq\dP$ is generic over\/ $\rL$. 
Then the sentence\/ $\vpi[\xx G]$ is true in\/ $\rL[G]$  
if and only if\/ $\sus T\in G\,(T\fo\vpi)$.
\ele
\bpf
{\ubf Base of induction:} 
$\vpi$ is arithmetic or belongs to $\lis11\cup\lip11$,
as in \ref{d:fo}\ref{fo1}. 
If $T\in G$ and $T\fo\vpi$, then $\vpi[\xx G]$ 
holds by the Shoenfield absoluteness theorem, as 
$\xx G\in[T]$. 
The inverse holds by Lemma~\ref{dens}. 

{\ubf Step $\lip1n\imp\lis1{n+1}$.} 
Let $\vpi$ be $\sus x\,\psi(x)$ where $\psi$ 
is of type $\lip1n$. 
Assume that $T\in G$ and $T\fo\vpi$. 
Then by Definition~\ref{d:fo}\ref{fo2}  
there is a code $\rc\in\ccf\cap\rL$ such that 
$T\fo\psi(\rc)$.
By the inductive hypothesis, 
the formula $\psi(\rc)[\xx G]$, that is,
$\psi[G](f_\rc(\xx G))$, is true in $\rL[G]$. 
But then $\vpi[\xx G]$ is obviously true as well. 

Conversely assume that $\vpi[\xx G]$ is true. 
Then there is a real $y\in\rL[G]\cap\dn$ such that
$\psi[G](y)$ is true.
By Lemma~\ref{repd},  
$y=f_\rc(\xx G)$ for a code $\rc\in\ccf\cap\rL$. 
But then $\psi(\rc)[\xx G]$ is true in $\rL[G]$. 
By the inductive hypothesis, there is a tree $T\in G$  
satisfying $T\fo\psi(\rc)$. 
Then $T\fo\vpi$ as well. 

{\ubf Step $\lis1n\imp\lip1{n}$.} 
Let $\vpi$ be a $\lip1n$ formula, $n\ge2$. 
By Lemma~\ref{dens}, there is a tree $T\in G$  
such that either $T\fo\vpi$ or $T\fo\vpi^-$.
If $T\fo\vpi^-$ then $\vpi^-[\xx G]$ is true by  
the inductive hypothesis, hence $\vpi[\xx G]$ is false. 
Now assume that $T\fo\vpi$. 
We have to prove that $\vpi[\xx G]$ is true. 
Suppose otherwise. 
Then $\vpi^-[\xx G]$ is true. 
By the inductive hypothesis, there is a tree  
$S\in G$ such that $S\fo\vpi^-$.
But the trees $S,T$ belong to the same generic set 
$G$, hence they are compatible, which leads to a 
contradiction with the assumption $T\fo\vpi$, 
according to Definition~\ref{d:fo}\ref{fo3}.
\epf

\ble
\lam{lun}
If a set\/ 
$G\sq\dP$ is\/ \dd\dP generic over\/ $\rL$, 
then it is true in\/ $\rL[G]$, that every countable\/ 
$\is1\nn$ set\/ $X\sq\dn$ satisfies\/ $Y\in\rL$. 
\ele 

\bpf 
We work in the context of Definition~\ref{vdp}. 
The first part of the proof is to show that $Y\sq\rL$. 
Suppose to the contrary that it holds in $\rL[G]$ 
that $Y\sq\dn$ is a countable $\is1\nn$ set, 
but $Y\not\sq\rL$. 
We have $Y=\ens{y\in\dn}{\vpi(y)}$  where  
$\vpi(y):=\sus z\,\psi(y,z)$ 
is a $\is1\nn$ formula, that is, $\lis1\nn$ formula 
without codes in $\ccf$. 
%
There is a tree $T_0\in G$ which \dd\dP forces that  
``$\ens{y\in\dn}{\vpi(y)}$ is a countable set and  
$\sus y(\vpi(y)\land y\nin\rL)$''. 
Our goal is to derive a contradiction.

By Lemma \ref{repd}, 
there exist codes $\rc,\rd\in\ccf\cap\rL$ such that the real  
$y_0=\rc[G]=f_\rc(\xx G)$ belongs to $Y$, so that  
$y_0\nin\rL$ and  
$\vpi(\rc)[\xx G]$ holds, that is, $\sus z\,\psi(\rc,z)[\xx G]$, 
and finally $\rd$ witnesses the existence quantifier, so that 
the sentence $\psi(\rc,\rd)[\xx G]$ holds in $\rL[G]$. 
By Lemma~\ref{ft.} as $\psi$ is a $\lip1{\nn-1}$ formula, 
there is a tree $T_1\in G$ satisfying  
$T_1\fo\psi(\rc,\rd)$.

We can wlog assume that $T_1\sq T_0$ and that,  
in $\rL$, the map $f_\rc$ is either a constant or a 
bijection on $[T_1]$, by Lemma~\ref{jden}\ref{prop4}. 

{\ubf Case 1:} $f_\rc\res[T_1]$ is a constant, that is,  
there exists a real $y_1\in\dn\cap\rL$ such that  
$f_\rc(x)=y_1$ for all $x\in[T_1]$. 
But then $y_1= f_\rc(\xx G)=y_0$, however $y_0\nin\rL$  
while $y_1\in\rL$, which is a contradiction.

{\ubf Case 2:} $f_\rc\res[T_1]$ is a bijection.
As $T_1$ is a Silver tree, the set 
$H=\oi{T_1}\sq\om$ of all its splitting levels 
is infinite (Definition~\ref{std}). 
Let $h\in\rL$, $h\sq H$. 
Then $h\app T_1=T_1$, and we have 
$T_1\fo \psi(h\app\rc,h\app\rd)$ 
By Lemma~\ref{inv}. 
Therefore the formula $\psi(h\app\rc,h\app\rd)[\xx G]$ 
is true in $\rL[G]$ by Lemma~\ref{ft.}. 
But this formula coincides with  
$\psi(f_{h\app\rc}(\xx G),f_{h\app\rd}(\xx D))$, 
hence we have 
$\vpi(f_{h\app\rc}(\xx G))$ in $\rL[G]$. 
This implies  
$f_{h\app\rc}(\xx G)\in Y$, or equivalently, 
$f_\rc(h\app\xx G)\in Y$. 

However if sets $h,h'\in\rL$, $h\cup h'\sq H$,
satisfy $h\ne h'$ then $h\app\xx G\ne h'\app\xx G$, 
and hence $f_\rc(h\app\xx G)\ne f_\rc(h'\app\xx G)$, 
as $f_\rc$ is a bijection on $[T_1]$ 
(and $\xx G\in[T_1]$).  
Thus picking different sets $h\in\rL\cap\pws H$ we get  
uncountably many different elements of the set $Y$ in  
$\rL[G]$, which contradicts to the choice of $Y$.

The proof of $Y\sq\rL$ is accomplished. 

To prove $Y\in\rL$, a stronger statement, it suffices 
now to show that if $y_0\in\dn\cap\rL$ then 
$y_0\in Y$ iff $\sus T\in\pes\,(T\fo\vpi(\rc_0))$, where 
$\rc_0\in\ccf\cap\rL$ is the code of the constant 
function $f_{\rc_0}(x)=y_0$, $\kaz x\in\dn$.

If $y_0\in Y$ then the formula $\vpi(y_0)$, equal 
to $\vpi(\rc_0)[x]$ for any $x$, is true in $\rL[G]$ 
by the choice of $\vpi$. 
It follows by Lemma~\ref{ft.} that there is a tree 
$T\in G$ satisfying $T\fo\vpi(\rc_0)$, as required. 

Now suppose that $T\in\pes$ 
(not necessarily $\in\dP$!) 
and $T\fo\vpi(\rc_0)$. 
As the set 
$D=\ens{T\in\pes}{\oi T\text{ is co-infinite}}$ 
(see Lemma~\ref{jden}\ref{prop5})
is open dense in $\pes$, 
we can assume that $\oi T$ is co-infinite. 
On the other hand, it follows from 
Lemma~\ref{jden}\ref{prop5} that there is a 
tree $S\in G\cap D$, 
so that $\oi S$ is co-infinite as well. 
Now we have $S\fo\vpi(\rc_0)$ by Corollary~\ref{tt'}, 
and then $\vpi(\rc_0)[\xx G]$ is true in $\rL[G]$ 
by Lemma~\ref{ft.}, that is, $\vpi(y_0)$ holds 
in $\rL[G]$, and $y_0\in Y$, as required.  
\epf

\bpf[\ubf Theorem~\ref{Tun}, the main theorem]
We assert that any \dd\dP generic extension 
$\rL[G]=\rL[\xx G]$ 
satisfies conditions \ref{Tun1}, \ref{Tun2}, \ref{Tun3} 
of the theorem.
That $\xx G\nin\od$ in \ref{Tun1} follows 
 by Lemma~\ref{sym}. 
The minimality follows from Lemma~\ref{jden}\ref{prop4}
by Lemma~\ref{repd} (continuous reading of names).
We further have \ref{Tun2} by Corollary~\ref{gdef}, 
and we have \ref{Tun3} by Lemma~\ref{lun}.
\epf


\bibliographystyle{plain}
\addcontentsline{toc}{subsection}{\hspace*{5.5ex}References}
{\small

\bibliography{34}
}

\def\indexname{\large Index%
\addcontentsline{toc}{subsection}{\hspace*{5.5ex}Index}}
\small\printindex

\end{document}